\documentclass[12pt]{article}

\usepackage{hyperref}
\usepackage{amsmath}
\usepackage{amssymb}
\usepackage{amsfonts}
\usepackage{amsopn}
\usepackage{amsthm}
\usepackage[all]{xy}

\swapnumbers
\theoremstyle{plain}
\newtheorem{thm}{Theorem}[section]
\newtheorem*{mainthm}{Theorem}
\newtheorem{lem}[thm]{Lemma}
\newtheorem{cor}[thm]{Corollary}
\newtheorem{prop}[thm]{Proposition}

\theoremstyle{definition}
\newtheorem{ntt}[thm]{}
\newtheorem{ex}[thm]{Example}
\newtheorem{rem}[thm]{Remark}
\newtheorem{dfn}[thm]{Definition}

\newcommand{\xra}[1]{\xrightarrow{#1}}   
\newcommand{\zz}{\mathbb{Z}}   
\newcommand{\qq}{\mathbb{Q}}   

\newcommand{\A}{\mathrm{A}}    
\newcommand{\B}{\mathrm{B}}    
\newcommand{\C}{\mathrm{C}}    
\newcommand{\D}{\mathrm{D}}    
\newcommand{\E}{\mathrm{E}}    
\newcommand{\F}{\mathrm{F}_4}  
\newcommand{\G}{\mathrm{G}_2}  

\DeclareMathOperator{\SO}{\mathrm{SO}}  
\DeclareMathOperator{\SL}{\mathrm{SL}}      
\DeclareMathOperator{\PGO}{\mathrm{PGO}}  
\DeclareMathOperator{\OO}{\mathrm{O}^+}  
\DeclareMathOperator{\PGSp}{\mathrm{PGSp}}  
\DeclareMathOperator{\Spin}{\mathrm{Spin}}  

\newcommand{\id}{\mathrm{id}}       
\newcommand{\pr}{\mathrm{pr}}       
\newcommand{\res}{\mathrm{res}}     
\newcommand{\cores}{\mathrm{cores}} 
\newcommand{\BX}{\overline{X}}   

\newcommand{\BG}{\overline{G}}   

\DeclareMathOperator{\im}{\mathrm{im}}      
\DeclareMathOperator{\End}{\mathrm{End}}    
\DeclareMathOperator{\SB}{\mathrm{SB}}      

\DeclareMathOperator{\M}{\mathcal{M}}     
\DeclareMathOperator{\RR}{\mathcal{R}}    
\DeclareMathOperator{\rr}{\mathfrak{r}}    
\DeclareMathOperator{\Spec}{\mathrm{Spec}}  

\DeclareMathOperator{\CH}{\mathrm{CH}}      
\DeclareMathOperator{\Ch}{\mathrm{Ch}}      
\DeclareMathOperator{\CHO}{\overline{\CH}}  
\DeclareMathOperator{\ChO}{\overline{\Ch}}  
\newcommand{\ChGi}{\mathop{\mathfrak{Ch}}(G)} 

\newcommand{\XX}{\mathfrak{X}}  
\newcommand{\BXX}{\overline{\XX}}

\DeclareMathOperator{\HH}{\mathrm{H}}      
\DeclareMathOperator{\rk}{\mathrm{rk}} 
\DeclareMathOperator{\tr}{\mathrm{tr}} 
\DeclareMathOperator{\ind}{\mathrm{ind}} 
\DeclareMathOperator{\codim}{\mathrm{codim}}  

\title{$J$-invariant of linear algebraic groups}

\author{V.~Petrov, N.~Semenov, K.~Zainoulline\footnote{The paper is based on the PhD Thesis of the first author.
Partially supported by CNRS, DAAD A/04/00348, INTAS, SFB701, DFG GI706/1-1.}}

\date{}

\begin{document}

\maketitle

\begin{abstract}
Let $G$ be a semisimple linear algebraic group of inner type 
over a field $F$ and $X$ be a 
projective homogeneous $G$-variety such that $G$ splits
over the function field of $X$.
In the present paper we introduce an invariant of $G$ called
$J$-invariant 
which characterizes the motivic behavior of $X$.
This generalizes the respective notion invented by A.~Vishik in the
context of quadratic forms.
As a main application
we obtain a uniform proof of
all known motivic decompositions of generically split
projective homogeneous varieties 
(Severi-Brauer varieties, Pfister quadrics, maximal orthogonal Grassmannians,
$\G$- and $\F$-varieties)
as well as provide new examples (exceptional varieties of types
$\E_6$, $\E_7$ and $\E_8$).
We also discuss relations with
torsion indices, canonical dimensions
and cohomological invariants of the group $G$.  
\end{abstract}

\section*{Introduction}
Let $G$ be a semisimple linear algebraic group over a field $F$
and $X$ be a projective homogeneous $G$-variety.
In the present paper we address the problem of computing 
the Grothendieck-Chow motive $\M(X)$ of $X$ or, in other words,
providing a direct sum decomposition of $\M(X)$.

This problem turns out to be strongly related with
several classical  conjectures concerning algebraic cycles.
For instance, the motivic decomposition of a Pfister quadric
plays a major role in the proof of Milnor's conjecture by V.~Voevodsky.
The proof of the generalization of this conjecture 
known as the Bloch-Kato conjecture was recently announced 
by M.~Rost and V.~Voevodsky. It essentially uses 
motivic decompositions of the norm varieties which 
are closely related to projective homogeneous varieties. 

Another deep application deals with the famous Kaplansky problem on
the values of the $u$-invariant of a field. 
It has a long history starting
from the works of A.~Merkurjev and O.~Izhboldin.
Recently an essential breakthrough 
in this problem was achieved by
A.~Vishik \cite{Vi06}, where he used the $J$-invariant
of an orthogonal group.
The present paper was mostly motivated by this result.
The invariant that we introduce and study is
a generalization of the $J$-invariant of A.~Vishik to
an arbitrary semisimple algebraic group.

It was first observed by B.~K\"ock \cite{Ko91} that
if the group $G$ is split, i.e.,
contains a split maximal torus, then 
the motive of $X$ has the simplest possible
decomposition -- it is isomorphic to
a direct sum of twisted Tate motives.
The next step was done by 
V.~Chernousov, S.~Gille and A.~Merkurjev \cite{CGM} and P.~Brosnan 
\cite{Br05}. 
They proved that 
if $G$ is isotropic, i.e.,
contains a split 1-dimensional torus, then the motive of $X$ 
can always be decomposed as a direct sum of
the motives of projective homogeneous varieties of smaller dimensions 
corresponding
to anisotropic groups, thus,
reducing the problem to the anisotropic case.

For anisotropic groups only very few partial results are known.
In this case the components 
of a motivic decomposition of $X$ 
are expected to have a non-geometric nature, i.e., can not be identified
with (twisted) motives of some other varieties. 
The first examples of such decompositions were provided
by M.~Rost \cite{Ro98}. He proved that the motive
of a Pfister quadric decomposes as a direct sum
of twisted copies of a certain a priori non geometric
motive $\RR$ called Rost motive.
The motives of Severi-Brauer varieties were computed by N.~Karpenko
\cite{Ka96}.
For exceptional varieties examples of motivic decompositions 
were provided by J.-P. Bonnet \cite{Bo03} 
(varieties of type $\G$) and 
by S.~Nikolenko, N.~Semenov, K.~Zainoulline \cite{NSZ} 
(varieties of type $\F$). Observe that in all these examples
the respective group $G$ splits over the generic point of $X$.
Such varieties will be called {\em generically split}.

In the present paper we provide a uniform proof of all these results.
Namely, we prove that (see Theorem~\ref{thmMain})

\begin{mainthm} 
Let $G$ be a semisimple linear algebraic group 
of inner type over
a field $F$ and $p$ be a prime integer. 
Let $X$ be a generically split 
projective homogeneous $G$-variety.
Then the Chow motive of $X$ with $\zz/p$-coefficients
is isomorphic to a direct sum
$$
\M(X;\zz/p)\simeq \bigoplus_{i\in\mathcal{I}}\RR_p(G)(i)
$$
of twisted copies of an indecomposable motive $\RR_p(G)$ for some
finite multiset $\mathcal{I}$ of non-negative integers. 
\end{mainthm}

Observe that the motive $\RR_p(G)$ depends only on 
$G$ and $p$ but
not on the type of a parabolic subgroup defining $X$.
Moreover, considered with $\qq$-coefficients 
it always splits as a direct sum of twisted Tate motives. 

Our proof is based on two different observations.
The first is the {\it Rost Nilpotence Theorem}. 
It was originally proven for projective 
quadrics by M.~Rost and then generalized to arbitrary 
projective homogeneous varieties by 
P.~Brosnan \cite{Br05}, V.~Chernousov, S.~Gille and A.~Merkurjev \cite{CGM}.
Roughly speaking, this result plays a role of the Galois descent
for motivic decompositions over a separable closure $\bar F$ of $F$.
Namely, it reduces the problem to the 
{\em description of idempotent cycles} in the endomorphism group
$\End (\M(X_{\bar F};\zz/p))$ 
which are defined over $F$. 

To provide such cycles we use the second observation 
which comes from the topology of compact Lie groups. 
In paper \cite{Kc85} V.~Kac 
invented the notion of {\em $p$-exceptional degrees} --
the numbers which relate the degrees of mod $p$ basic polynomial invariants
and the $p$-torsion part of the Chow ring of a compact Lie group.
These numbers have combinatorial nature.
By the result of K.~Zainoulline \cite{Za06}
there is a strong interrelation between 
$p$-exceptional degrees and the subgroup
of cycles in $\End(\M(X_{\bar F};\zz/p))$ defined over $F$.
To describe this subgroup we introduce the notion
of the {\em $J$-invariant}
of a group $G$ mod $p$ denoted by $J_p(G)$
(see Definition~\ref{jinv}).
In the most cases the values
of $J_p(G)$ were implicitly computed by V.~Kac in \cite{Kc85} and can
easily be extracted from Table~\ref{bigtab}. 

It follows from the proof that
the $J$-invariant measures the `size' of
the motive $\RR_p(G)$ 
and, hence, characterizes
the motivic decomposition of $X$. 
Observe that if the $J$-invariant takes its minimal possible non-trivial
value $J_p(G)=(1)$, then the motive $\RR_p(G)\otimes\qq$ has the following
recognizable decomposition
(cf. \cite[\S5]{Vo03} and \cite[\S5]{Ro06})
$$
\RR_p(G)\otimes \qq \simeq \bigoplus_{i=0}^{p-1} 
\qq(i\cdot \tfrac{p^{n-1}-1}{p-1}),\text{ where }n=2\text{ or }3. 
$$
The assignment $G\mapsto\RR_p(G)$ can be viewed as a
motivic analog of the {\em cohomological invariant} of $G$
given by the Tits class of $G$ if $n=2$ and
by the Rost invariant of $G$ if $n=3$. In these cases the motive $\RR_p(G)$
coincides with a {\em generalized Rost motive}.

We also generalize some of the results of paper \cite{CPSZ}.
Namely, using the motivic version of the result
of D.~Edidin and W.~Graham \cite{EG97} 
on {\em cellular fibrations} we provide
a general formula which
expresses the motive of the total space of a cellular fibration 
in terms of the motives of its base (see Theorem~\ref{thmf}). 
We also provide several criteria for the existence of {\em liftings
of motivic decompositions} via the reduction map $\zz\to \zz/m$.
We prove that such liftings always exist (see Theorem~\ref{lem47}).

The paper is organized as follows.
In the first section we provide several auxiliary
facts concerning motives and rational cycles.
Rather technical section~\ref{IDLI} is devoted to lifting of idempotents.
In section~\ref{EDGR} we discuss the motives of cellular fibrations.
The next section is devoted to the notion of a $J$-invariant.
The proof of the main result is given in 
section~\ref{MDEC}.
The last two sections are devoted 
to various applications of the $J$-invariant
and examples of motivic decompositions.
In particular, we discuss the relations with canonical $p$-dimensions,
degrees of zero-cycles and the Rost invariant. 


\section{Chow motives and rational cycles}

In the present section we provide several auxiliary facts
concerning algebraic cycles, correspondences and motives 
which will be extensively used in the sequel. 
We follow the notation and definitions from \cite[Ch.~XII]{EKM}
(see also \cite{Ma68}).

Let $X$ be a smooth projective irreducible variety over a field $F$. 
Let $\CH_i(X;\Lambda)$ be the Chow group
of cycles of dimension $i$ on $X$ 
with coefficients in a commutative ring $\Lambda$.
For simplicity we denote $\CH(X;\zz)$ by $\CH(X)$.

\begin{ntt} \label{chmot}
Following \cite[\S63]{EKM}
an element $\phi\in \CH_{\dim X+d}(X\times Y;\Lambda)$ 
is called a {\em correspondence} between $X$ and $Y$ of degree $d$ with
coefficients in $\Lambda$. 
Let $\phi \in\CH(X\times Y;\Lambda)$ and 
$\psi\in \CH(Y\times Z;\Lambda)$ be correspondences
of degrees $d$ and $e$ respectively.
Then their product $\psi\circ\phi$ is defined by the formula 
$(\pr_{XZ})_*(\pr_{XY}^*(\phi)\cdot\pr_{YZ}^*(\psi))$
and has degree $d+e$.
The {\em correspondence product} endows the group 
$\CH(X\times X;\Lambda)$ with a ring structure. 
The identity element of this ring
is the class of the diagonal $\Delta_X$.
Given $\phi\in\CH(X\times X;\Lambda)$ of degree $d$ 
we define a $\Lambda$-linear
map $\CH_i(X;\Lambda)\to \CH_{i+d}(Y;\Lambda)$ by
$\alpha\mapsto (\pr_Y)_*(\pr_X^*(\alpha)\cdot \phi)$.
This map is called {\em realization} of $\phi$ and is denoted
by $\phi_\star$.
By definition
$(\psi\circ\phi)_\star=\psi_\star\circ\phi_\star$.
Given a correspondence $\phi$ we denote by $\phi^t$ its {\em transpose}.
\end{ntt}

\begin{ntt}\label{chmm} 
Following \cite[\S~64]{EKM}
let $\M(X;\Lambda)$ 
denote the {\em Chow motive} of $X$ with
$\Lambda$-coefficients
and $\M(X;\Lambda)(n)=\M(X;\Lambda)\otimes \Lambda(n)$ denote the respective 
{\em twist} by the Tate motive.
For simplicity by $\M(X)$
we will denote $\M(X;\zz)$.
Recall that morphisms between
$\M(X;\Lambda)(n)$ and $\M(Y;\Lambda)(m)$
are given by correspondences 
of degree $n-m$ between $X$ and $Y$.
Hence, the group of endomorphisms
$\End(\M(X;\Lambda))$ coincides 
with the Chow group 
$\CH_{\dim X}(X\times X;\Lambda)$. 
Observe that 
to provide a direct sum decomposition of $\M(X;\Lambda)$ is the same
as to provide a family of pair-wise orthogonal idempotents 
$\phi_i\in\End(\M(X;\Lambda))$ such that $\sum_i \phi_i=\Delta_X$.
\end{ntt}

\begin{ntt} Assume that a motive $M$
is a direct sum of twisted Tate motives. In this case its
Chow group $\CH(M)$ is a free abelian group.
We define its {\em Poincar\'e polynomial} as
$$
P(M,t)=\sum_{i\ge 0} a_i t^i, 
$$
where $a_i$ is the rank of $\CH_i(M)$.
\end{ntt}

\begin{dfn}\label{spvar} 
Let $L/F$ be a field extension.
We say $L$ is a {\em splitting field} of a smooth projective variety $X$ 
or, equivalently, a variety $X$ {\em splits} over $L$ if 
the motive $\M(X;\zz)$ splits over $L$ as 
a finite direct sum of twisted Tate motives.
\end{dfn}

\begin{ex}\label{exspl}
A variety $X$ over a field $F$ is called {\em cellular} if
$X$ has a proper descending filtration by
closed subvarieties $X_i$ such that each complement $X_i\setminus X_{i+1}$
is a disjoint union of affine spaces defined over $F$. 
According to \cite[Corollary~66.4]{EKM} if $X$ is cellular,
then $X$ splits over $F$.

In particular, 
let $G$ be a semisimple linear algebraic group over a field $F$
and $X$ be a projective homogeneous $G$-variety.
Assume that the group $G$
splits over the generic point of $X$, i.e., $G_{F(X)}=G\times_F F(X)$
contains a split maximal torus defined over $F(X)$. 
Then $X_{F(X)}$ is a cellular variety
and, therefore, $F(X)$ is a splitting field of $X$. 
Some concrete examples of such varieties are provided in \ref{intmot}.
\end{ex}

\begin{ntt}
Assume $X$ has a splitting field $L$. We will write
$\CH(\overline X;\Lambda)$ for $\CH(X_L;\Lambda)$
and $\CHO(X;\Lambda)$ for the image of the restriction map
$\CH(X;\Lambda) \to \CH(\overline X;\Lambda)$  
(cf. \cite[1.2]{KM05}). 
Similarly, by $\M(\overline X;\Lambda)$ we denote the motive 
of $X$ considered over $L$.
If $M$ is a direct summand of $\M(X;\Lambda)$,
by $\overline M$ we denote the motive $M_L$.
The elements of $\CHO(X;\Lambda)$ 
will be called {\em rational} cycles on $X_L$ with respect
to the field extension $L/F$ and the coefficient ring $\Lambda$. 
If $L'$ is another splitting field of $X$, then
there is a chain of canonical isomorphisms
$\CH(X_L)\simeq \CH(X_{LL'})\simeq \CH(X_{L'})$, where $LL'$ is the composite
of $L$ and $L'$. Hence, the groups $\CH(\overline X)$ and $\CHO(X)$ do not
depend on the choice of $L$.
\end{ntt}

\begin{ntt}\label{remsp} According to \cite[Remark~5.6]{KM05} 
there is the K\"unneth decomposition 
$\CH(\overline X\times \overline X)=\CH(\overline X)\otimes \CH(\overline X)$ 
and Poincar\'e duality holds for $\CH(\overline X)$. 
The latter means that given a
basis of $\CH(\overline X)$ 
there is a dual one with respect to the non-degenerate pairing  
$(\alpha,\beta)\mapsto \deg(\alpha\cdot\beta)$, where $\deg$
is the degree map.
In view of the K\"unneth decomposition 
the correspondence product of cycles in $\CH(\overline X\times\overline X)$
is given by the formula
$(\alpha_1\times \beta_1)\circ (\alpha_2\times\beta_2)=
\deg(\alpha_1\beta_2)(\alpha_2\times\beta_1)$, the realization by
$(\alpha\times\beta)_\star(\gamma)=\deg(\alpha\gamma)\beta$ 
and the transpose by
$(\alpha\times\beta)^t=\beta\times\alpha$.
\end{ntt}

Sometimes we will use contravariant notation $\CH^*$ 
for Chow groups meaning
$\CH^i(X)=\CH_{\dim X-i}(X)$ for irreducible $X$.
The following important fact will be used in the proof of
the main theorem (see Lemma~\ref{sigrat}).

\begin{lem}\label{gensp}
Let $X$ and $Y$ be two smooth projective varieties such that 
$Y$ is irreducible, $F(Y)$ is a splitting field of $X$ and 
$Y$ has a splitting field.
For any $r$ consider the projection in the K\"unneth decomposition
$$
\pr_{0}\colon 
\CH^r(\overline X\times \overline Y)=
\bigoplus_{i=0}^r \CH^{r-i}(\overline X)\otimes \CH^i(\overline Y)
\to \CH^r(\overline X).
$$
Then for any $\rho \in \CH^r(\overline X)$ we have
$\pr_{0}^{-1}(\rho)\cap \overline{\CH^r}(X\times Y)\neq \emptyset$. 
\end{lem}

\begin{proof} Let $L$ be a common splitting field of $X$ and $Y$.
The lemma follows from the commutative diagram
$$
\xymatrix{
\CH^r(X\times_F Y) \ar[r]^{\res_{L/F}}\ar@{>>}[d]& 
\CH^r(X_L\times_L Y_L)\ar@{>>}[d] 
 \ar[dr]^{\pr_{0}}& \\
\CH^r(X_{F(Y)})\ar[r]^-\simeq & \CH^r((X_L)_{L(Y_L)}) 
 & \CH^r(X_L)\ar[l]_-\simeq
}
$$
where the left square is obtained 
by taking the generic fiber of the base change morphism $X_L\to X$;
the vertical arrows are taken from 
the localization sequence for Chow groups
and, hence, are surjective; and
the bottom horizontal maps are isomorphisms since $L$ is a splitting field.
\end{proof}

\begin{dfn}\label{rpres}
We say that a field extension $E/F$ is {\em rank preserving} with respect to 
$X$ if the restriction map $\res_{E/F}\colon\CH(X)\to\CH(X_E)$ 
becomes an isomorphism after tensoring with $\mathbb{Q}$.
\end{dfn}

\begin{lem}\label{prat} 
Assume $X$ has a splitting field. Then for any rank preserving finite field 
extension $E/F$ we have
$[E:F]\cdot\CHO(X_E)\subset\CHO(X)$.
\end{lem}
\begin{proof} 
Let $L$ be a splitting field containing $E$.  Let $\gamma$ be any element in 
$\CHO(X_E)$. By definition there exists $\alpha\in\CH(X_E)$ 
such that  $\gamma=\res_{L/E}(\alpha)$.
Since $\res_{E/F}\otimes \mathbb{Q}$ is an isomorphism, there exists an 
element $\beta\in\CH(X)$ and  a non-zero integer $n$ such that $\res_{E/F}
(\beta)=n\alpha$. By the projection formula 
$$
n\cdot \cores_{E/F}(\alpha)=\cores_{E/F}(\res_{E/F}(\beta))=[E:F]\cdot \beta.
$$
Applying $\res_{L/E}$ to the both sides of the identity we obtain
$$n(\res_{L/E}(\cores_{E/F}(\alpha)))=n[E:F]\cdot\gamma.$$ 
Therefore, $\res_{L/E}(\cores_{E/F}(\alpha))=[E:F]\cdot\gamma$. 
\end{proof}

We provide now examples of varieties for which any field extension 
is rank preserving and, hence, Lemma~\ref{prat} holds.

\begin{ntt}\label{nota}
Let $G$ be a semisimple linear algebraic group over a field $F$, 
$X$ be a projective homogeneous $G$-variety. Denote by $\mathcal D$ the 
Dynkin diagram of $G$. According to \cite{Ti66}
one can always choose a quasi-split 
group $G_0$ over $F$ with the same Dynkin diagram, 
a parabolic subgroup $P$ of $G_0$ and 
a cocycle $\xi\in H^1(F,G_0)$ such that 
$G$ is isogenic to the twisted form ${}_\xi G_0$ and 
$X$ is isomorphic to ${}_\xi(G_0/P)$.
If $G_0$ is split (see \ref{exspl}), then $G$ 
is called a group of \emph{inner type} over $F$.
\end{ntt}

\begin{lem}\label{surj} 
Let $G$ be a semisimple linear algebraic 
group of inner type over a field $F$
and $X$ a projective homogeneous $G$-variety. 
Then any field extension $E/F$ is rank preserving with respect to $X$ and 
$X\times X$.
\end{lem}

\begin{proof}
By  \cite[Theorem~2.2 and 4.2]{Pa94} the restriction map
$K_0(X) \to K_0(X_E)$  becomes an isomorphism after tensoring with 
$\qq$. Now the Chern character
$ch\colon K_0(X)\otimes\qq\to \CH^*(X)\otimes\qq$
is an isomorphism and respects pull-backs, hence, $E$ is rank preserving with 
respect to $X$. It remains to note that $X\times X$ is a homogeneous $G\times 
G$-variety.
\end{proof}

\begin{rem} For even dimensional quadrics with non-trivial
discriminant the restriction map $K_0(X)\to K_0(\BX)$ is not surjective and
Lemma~\ref{surj} doesn't hold. 
\end{rem}


\section{Lifting of idempotents}\label{IDLI}

This section is devoted to lifting of idempotents 
and isomorphisms.
First, we treat the case of general graded algebras.
The main results here are Lemma~\ref{izvrat} and
Proposition~\ref{nilpcor}.
Then, assuming {\em Rost Nilpotence} \ref{dfrnil}
we provide conditions
to lift motivic decompositions and isomorphisms (Theorem~\ref{lem47}).

\begin{ntt}
Let $A^*$ be a $\zz$-graded ring.
Assume we are given two orthogonal idempotents $\phi_i$ and $\phi_j$ in 
$A^0$ that is $\phi_i\phi_j=\phi_j\phi_i=0$.
We say an element $\theta_{ij}$  provides 
an {\em isomorphism} of degree $d$ between
idempotents $\phi_i$ and $\phi_j$  if 
$\theta_{ij}\in \phi_j A^{-d} \phi_i$ and 
there exists $\theta_{ji}\in \phi_i A^d \phi_j$  such that
$\theta_{ij} \theta_{ji}=\phi_j$ and  $\theta_{ji} \theta_{ij}=\phi_i$.
\end{ntt}

\begin{ex}\label{motex} 
Let $X$ be a smooth projective irreducible variety over a field $F$ 
and $\CH^*(X\times X;\Lambda)$ 
the Chow ring with coefficients in a commutative ring $\Lambda$.
Set $A^*=\End^*(\M(X;\Lambda))$, where
$$\End^{-i}(\M(X;\Lambda))=\CH^{\dim X-i}(X\times X;\Lambda)=
\CH_{\dim X+i}(X\times X;\Lambda),\quad i\in\zz $$ 
and
the multiplication is given by the correspondence product.
By definition $\End^0(\M(X;\Lambda))$ is the ring of endomorphisms
of the motive $\M(X;\Lambda)$ (see \ref{chmm}).
Note that a direct summand of $\M(X;\Lambda)$ can be identified
with a pair $(X,\phi_i)$, where
$\phi_i$ is an idempotent (see \cite[ch.~XII]{EKM}).
Then an isomorphism $\theta_{ij}$ of degree $d$ 
between $\phi_i$ and $\phi_j$ can be identified with 
an isomorphism between the motives $(X,\phi_i)$ and $(X,\phi_j)(d)$.
\end{ex}

\begin{dfn}\label{pres} 
Let $f\colon A^*\to B^*$ be a homomorphism of $\zz$-graded rings. 
We say that $f$ \emph{lifts decompositions} 
if given a family $\phi_i\in B^0$ of pair-wise orthogonal idempotents 
such that $\sum_i \phi_i=1_B$, 
there exists a family of pair-wise orthogonal idempotents $\varphi_i\in A^0$ 
such that  $\sum_i \varphi_i=1_A$ and 
each $f(\varphi_i)$ is isomorphic to $\phi_i$ 
by means of an isomorphism of degree $0$. 
We say $f$ \emph{lifts decompositions strictly} if, 
moreover, one can choose $\varphi_i$ such that 
$f(\varphi_i)=\phi_i$.
 
We say $f$ \emph{lifts isomorphisms}
if for any idempotents $\varphi_1$ and $\varphi_2$ in $A^0$ and 
any isomorphism $\theta_{12}$ of degree $d$ 
between idempotents $f(\varphi_1)$ and $f(\varphi_2)$ in $B^0$  
there exists an isomorphism $\vartheta_{12}$ of degree $d$ 
between $\varphi_1$ and $\varphi_2$. 
We say $f$ \emph{lifts isomorphisms strictly} if, 
moreover, one can choose $\vartheta_{12}$ such that 
$f(\vartheta_{12})=\theta_{12}$.
\end{dfn}

\begin{ntt}\label{prope}
By definition we have the following properties of morphisms
which lift decompositions and isomorphisms (strictly):
\begin{enumerate}
\item[(i)] 
Let $f\colon A^*\to B^*$ and 
$g\colon B^*\to C^*$ be two morphisms.
If both $f$ and $g$ lift decompositions or isomorphisms (strictly),
then so does the composite $g\circ f$.
\item[(ii)]
If $g\circ f$ lifts decompositions (resp. isomorphisms) 
and $g$ lifts isomorphisms, then $f$ lifts decompositions 
(resp. isomorphisms).
\item[(iii)]
Assume we are given a commutative diagram with $\ker f'\subset\im h$
$$
\xymatrix{
A^*\ar@{>>}[r]^f\ar@{^{(}->}[d]_h & B^*\ar@{^{(}->}[d]^{h'}\\
A'^*\ar@{>>}[r]^{f'} & B'^*.
}
$$
If $f'$ lifts {\em decompositions strictly} 
(resp. {\em isomorphisms strictly}), 
then so does $f$.
\end{enumerate}
\end{ntt}

\begin{lem}\label{izvrat}
Let $A$, $B$ be two rings, $A^0$, $B^0$ be their subrings,
$f^0\colon A^0\to B^0$ be a ring homomorphism, $f\colon A\to B$ be a map of 
sets satisfying the following conditions:
\begin{itemize}
\item $f(\alpha)f(\beta)$ equals either $f(\alpha\beta)$ or $0$ for all
$\alpha,\beta\in A$;
\item $f^0(\alpha)$ equals $f(\alpha)$ if $f(\alpha)\in B^0$ or $0$ otherwise;
\item $\ker f^0$ consists of nilpotent elements.
\end{itemize}
Let $\varphi_1$ and $\varphi_2$ be two idempotents in $A^0$, $\psi_{12}$ and 
$\psi_{21}$ be elements in $A$ such that $\psi_{12}A^0\psi_{21}\subset A^0$, 
$\psi_{21}A^0\psi_{12}\subset A^0$,
$f(\psi_{21})f(\psi_{12})=f(\varphi_1)$, 
$f(\psi_{12})f(\psi_{21})=f(\varphi_2)$. 

Then there exist elements
$\vartheta_{12}\in\varphi_2 A^0\psi_{12}A^0\varphi_1$ and
$\vartheta_{21}\in\varphi_1 A^0\psi_{21}A^0\varphi_2$ such that
$\vartheta_{21}\vartheta_{12}=\varphi_1$, 
$\vartheta_{12}\vartheta_{21}=\varphi_2$,
$f(\vartheta_{12})=f(\varphi_2)f(\psi_{12})=f(\psi_{12})f(\varphi_1)$, 
$f(\vartheta_{21})=f(\varphi_1)f(\psi_{21})=f(\psi_{21})f(\varphi_2)$.
\end{lem}
\begin{proof}
Since $\ker f^0$ consists of nilpotents, $f^0$ sends non-zero idempotents in 
$A^0$ to non-zero idempotents in $B^0$; in particular, $f(\varphi_1)=f^0
(\varphi_1)\ne 0$,
$f(\varphi_2)=f^0(\varphi_2)\ne 0$.
Observe that
$$
f(\psi_{12})f(\varphi_1)=f(\psi_{12})f(\psi_{21})f(\psi_{12})=f(\varphi_2)f(\psi_{12})
$$
and, similarly, $f(\psi_{21})f(\varphi_2)=f(\varphi_1)f(\psi_{21})$.
Changing $\psi_{12}$ to $\varphi_2\psi_{12}\varphi_1$ and $\psi_{21}$ to 
$\varphi_1\psi_{21}\varphi_2$ we may assume that 
$\psi_{12}\in\varphi_2A\varphi_1$ and $\psi_{21}\in\varphi_1A\varphi_2$. We 
have
$$
f^0(\varphi_2)=f(\varphi_2)=f(\psi_{12})f(\psi_{21})=f(\psi_{12}\psi_{21})
=f^0(\psi_{12}\psi_{21}).$$
Therefore $\alpha=\psi_{12}\psi_{21}-\varphi_2\in A^0$ is nilpotent, say 
$\alpha^n=0$. Note that $\varphi_2\alpha=\alpha=\alpha\varphi_2$. Set 
$\alpha^\vee=\varphi_2-\alpha+\ldots+(-1)^{n-1}\alpha^{n-1}\in A^0$; then 
$\alpha\alpha^\vee=\varphi_2-\alpha^\vee$, 
$\varphi_2\alpha^\vee=\alpha^\vee=\alpha^\vee\varphi_2$ and $f(\varphi_2)=f^0
(\varphi_2)=f^0(\alpha^\vee)=f(\alpha^\vee)$. Therefore setting 
$\vartheta_{21}=\psi_{21}\alpha^\vee$ we have $\vartheta_{21}\in 
\varphi_1A\varphi_2$, 
$$
\psi_{12}\vartheta_{21}=\psi_{12}\psi_{21}\alpha^\vee=(\varphi_2+\alpha)\alpha^\vee=\alpha^\vee+\varphi_2-\alpha^\vee=\varphi_2$$
 and 
$f(\psi_{21})f(\varphi_2)=f(\psi_{21})f(\alpha^\vee)=f(\psi_{21}\alpha^\vee)=f(\vartheta_{21})$. This also implies that 
$\vartheta_{21}\psi_{12}$ is an idempotent. 

We have
$$
f^0(\varphi_1)=f(\varphi_1)=f(\vartheta_{21})f(\psi_{12})
=f(\vartheta_{21}\psi_{12})=f^0(\vartheta_{21}\psi_{12}),
$$
where the last equality holds, since 
$f(\vartheta_{21}\psi_{12})=f^0(\varphi_1)$ belongs to $B^0$ and 
$f^0$ satisfies the second condition.
Therefore $\beta=\vartheta_{21}\psi_{12}-\varphi_1\in A^0$ is nilpotent. Note 
that $\beta\varphi_1=\beta=\varphi_1\beta$. Now
$\varphi_1+\beta=(\varphi_1+\beta)^2=\varphi_1+2\beta+\beta^2$ and therefore
$\beta(1+\beta)=0$. But $1+\beta$ is invertible and hence we have $\beta=0$. 
It means that $\vartheta_{21}\psi_{12}=\varphi_1$ and we can set 
$\vartheta_{12}=\psi_{12}$.
\end{proof}

\begin{prop}\label{nilpcor}
Let $f\colon A^*\to B^*$ be a surjective homomorphism 
such that the kernel of the restriction of $f$ to $A^0$ 
consists of nilpotent elements. 
Then $f$ lifts decompositions and
isomorphisms strictly.
\end{prop}

\begin{proof}
The fact that $f$ lifts {\em decompositions strictly} follows
from \cite[Proposition~27.4]{AnFu}.

Let $\varphi_1$ and $\varphi_2$ be two idempotents in $A^0$
and $\theta_{12}$ be an isomorphism between $f(\varphi_1)$ and $f(\varphi_2)$.
Let $\psi_{12}$ in $A$ (resp. $\psi_{21}$)
be a homogeneous lifting of $\theta_{12}$ (resp. $\theta_{21}$).
The proposition follows now from Lemma~\ref{izvrat}.
\end{proof}

\begin{cor}\label{karpcor}
Let $m$ be an integer and $m=p_1^{n_1}\ldots p_l^{n_l}$ be its prime 
factorization.
Then the product of reduction maps 
$$\End^*(\M(X;\zz/m))\twoheadrightarrow \prod_{i=1}^l \End^*(\M(X;\zz/p_i))$$ 
lifts decompositions and isomorphisms strictly.
\end{cor}

\begin{proof} 
We apply Proposition~\ref{nilpcor} to the case
$A^*=\End^*(\M(X;\zz/p_i^{n_i}))$, $B^*=\End^*(\M(X;\zz/p_i))$
and
the reduction map $f_i\colon A^*\twoheadrightarrow B^*$. 
We obtain that $f_i$
lifts decompositions and isomorphisms strictly for each $i$.
To finish the proof observe that
by the Chinese remainder theorem 
$\End^*(\M(X;\zz/m))\simeq \prod_{i=1}^l \End^*(\M(X;\zz/p_i^{n_i}))$.
\end{proof}

We are coming to the following important definition.

\begin{dfn}\label{dfrnil} 
Let $X$ be a smooth projective variety over a field $F$. 
Assume that $X$ has a splitting field (see \ref{spvar}).
We say that {\em Rost Nilpotence} holds for $X$ if 
the kernel of the restriction map
$$
\res_E \colon\End^*(\M(X_E;\Lambda))\to\End^*(\M(\overline X;\Lambda))
$$
consists of nilpotent elements for all field 
extensions $E/F$ and all rings of coefficients $\Lambda$.
\end{dfn}

\begin{lem}\label{rostnil}
Let $X$ be a smooth projective variety which splits over any field 
over which it has a rational point. Then Rost Nilpotence holds for $X$.
\end{lem} 

\begin{proof} By \cite[Theorem~67.1]{EKM} if $\alpha$ is in the kernel
of the restriction map $\res_E$ then $\alpha^{\circ (\dim X+1)}=0$.
\end{proof}

\begin{lem}\label{rostcor}
Assume that Rost Nilpotence holds for $X$.
Then for any field extension $E/F$ the restriction 
$\res_E\colon \End^*(\M(X_E;\Lambda))\twoheadrightarrow\im(\res_E)$ 
onto the image 
lifts decompositions and isomorphisms strictly.
\end{lem}

\begin{proof}
Apply Proposition~\ref{nilpcor} to the homomorphism 
$\res_E\colon A^*\twoheadrightarrow B^*$
between the graded rings
$A^*=\End^*(\M(X_E;\Lambda))$ and
$B^*=\im(\res_E)$.
\end{proof}

\begin{cor}\label{transf} 
Assume that Rost Nilpotence holds for $X$. 
Let $m$ be an integer and 
$E/F$ be a field extension of degree coprime 
to $m$ which is rank preserving with respect to $X\times X$ 
(see \ref{rpres}). 
Then the restriction map 
$$\res_{E/F}\colon
\End^*(\M(X;\zz/m))\to\End^*(\M(X_E;\zz/m))$$ 
lifts decompositions and 
isomorphisms.
\end{cor}
\begin{proof}
By Lemma~\ref{prat} we have
$\im(\res_E)=\im(\res_F)$.
We apply now Lemma~\ref{rostcor} and \ref{prope}(ii) with
$A^*=\End^*(\M(X;\zz/m))$, $B^*=\End^*(\M(X_E;\zz/m))$ and
$C^*=\im(\res_E)$.
\end{proof}

\begin{dfn} Let $V^*$ be a free graded $\Lambda$-module of finite rank
and $A^*=\End^*(V^*)$ be its ring of endomorphisms, where
$\End^{-d}(V^*)$, $d\in \zz$, is the group of endomorphisms of $V^*$ 
decreasing the degree by $d$.
Assume we are given a direct sum decomposition of 
$V^*=\bigoplus_i \im \phi_i$ by means
of idempotents $\phi_i$ in $A^*$.
We say that this decomposition is $\Lambda$-{\em free} if
all graded components of $\im \phi_i$ are free $\Lambda$-modules.
Observe that if $\Lambda=\zz$ or $\zz/p$, where $p$ is prime,
then any decomposition is $\Lambda$-free.
\end{dfn}

\begin{ex}\label{spdcm} Assume $X$ has a splitting field.
Define $V^*=\CH^*(\overline X)$. 
Then by Poincar\'e duality and by the 
K\"unneth decomposition (see \ref{remsp}) we have
$\End^*(V^*)=\End^*(\M(\overline X))$ 
(see Example~\ref{motex}).
\end{ex}

\begin{lem}\label{surjlem} 
The map $\SL_l(\zz)\to \SL_l(\zz/m)$ induced by the reduction
modulo $m$ is surjective.
\end{lem}
\begin{proof} Since $\zz/m$ is a semi-local ring, 
the group $\SL_l(\zz/m)$ is generated by elementary matrices 
(see \cite[Theorem~4.3.9]{HOM89}).
\end{proof}

\begin{prop}\label{lem46}
Consider a free graded $\zz$-module $V^*$ of finite rank and
the reduction map $f\colon \End^*(V^*)\twoheadrightarrow 
\End^*(V^*\otimes_{\zz}\zz/m)$.
Then $f$ lifts $\zz/m$-free decompositions strictly.
Moreover, if $(\zz/m)^\times=\{\pm 1\}$, then 
$f$ lifts isomorphisms of $\zz/m$-free decompositions strictly.
\end{prop}

\begin{proof} 
We are given a decomposition $V^k\otimes_{\zz}\zz/m=\oplus_i W_i^k$,
where $W_i^k$ is the $k$-graded component of $\im\phi_i$.
Present $V^k$ as a direct sum  $V^k=\bigoplus_i V_i^k$ of free $\zz$-modules
such that $\rk_\zz V_i^k=\rk_{\zz/m} W_i^k$. Fix a $\zz$-basis
$\{v_{ij}^k\}_j$ of $V_i^k$. For each $W_i^k$ choose a basis
$\{w_{ij}^k\}_j$ such that the linear transformation $D^k$ of
$V^k\otimes_{\zz}\zz/m$ sending each $v_{ij}^k\otimes 1$ to $w_{ij}^k$ has
determinant $1$. By Lemma~\ref{surjlem} there is a lifting $\tilde D^k$ of
$D^k$ to a linear transformation of $V^k$. So we obtain
$V^k=\bigoplus_i \tilde W_i^k$, where  $\tilde W_i^k=\tilde D^k(V_i^k)$
satisfies $\tilde W_i^k\otimes_{\zz}\zz/m=W_i^k$. 
It remains to define $\varphi_i$ on each $V^k$ 
to be the projection onto $\tilde W_i^k$.

Now let $\varphi_1,\varphi_2$ be two idempotents in $\End^*(V^*)$. 
Denote by $V_i^k$ the $k$-graded component of $\im\varphi_i$. 
An isomorphism $\theta_{12}$ between $\varphi_1\otimes 1$ and $\varphi_2\otimes 1$ of degree $d$ can be identified with a family of isomorphisms $\theta_{12}^k\colon V_1^k\otimes\zz/m\to V_2^{k-d}\otimes\zz/m$. 
In the case $(\zz/m)^\times=\{\pm1\}$ all these isomorphisms are given
by matrices with determinants $\{\pm 1\}$ and, hence, can be 
lifted to isomorphisms $\vartheta_{12}^k\colon V_1^k\to V_2^{k-d}$
by Lemma~\ref{surjlem}.
\end{proof}

Now we are ready to state and to prove the main result of this section.

\begin{thm}\label{lem47} 
Let $X$ be a smooth projective irreducible variety over a field $F$.
Assume that $X$ has a splitting field of degree $m$ 
which is rank preserving with respect to $X\times X$. 
Assume that Rost Nilpotence holds for $X$. 
Consider only decompositions of $\M(X;\zz/m)$ which become $\zz/m$-free
over the splitting field.
Then the reduction map 
$$f\colon \End^*(\M(X;\zz))\twoheadrightarrow \End^*(\M(X;\zz/m))$$
lifts such decompositions.
If additionally $(\zz/m)^\times=\{\pm 1\}$, then this map 
lifts isomorphisms of such decompositions.
\end{thm}

\begin{proof}
Consider the diagram
$$
\xymatrix{
\End^*(\M(X;\zz))\ar@{>>}[r]\ar@{>>}[d]_-f & \im(\res_F)\ar@{^{(}->}[r]^-h\ar@{>>}[d]^-{\bar f} & 
\End^*(\M(\overline X;\zz))\ar@{>>}[d]^-{f'} \\
\End^*(\M(X;\zz/m)) \ar@{>>}[r] & \im(\res_F)\ar@{^{(}->}[r] & \End^*(\M(\overline X;\zz/m))
}
$$
Recall that
we can identify
$\End^{-d}(\M(\overline X))$ with the group of endomorphisms of 
$\CH^*(\overline X)$ which decrease the grading by $d$ (see \ref{spdcm}).
Applying Proposition~\ref{lem46} to the case $V^*=\CH^*(\overline X)$
we obtain that the map $f'$ lifts decompositions strictly. 
Moreover, if 
$(\zz/m)^\times=\{\pm 1\}$ then $f'$ lifts isomorphisms strictly.

By Lemma~\ref{prat} $\ker f' \subset \im h$ and, therefore, 
applying \ref{prope}(iii) we obtain that $\bar f$ lifts
decompositions strictly
and, moreover, $\bar f$ lifts isomorphisms strictly 
if $(\zz/m)^\times=\{\pm 1\}$.  

Now by Lemma~\ref{rostcor}
the horizontal arrows of the left square 
lift decompositions and isomorphisms strictly.
It remains to apply \ref{prope}(i) and (ii).
\end{proof}


\section{Motives of fibered spaces}\label{EDGR}

In the present section we discuss motives of cellular fibration.
The main result (Theorem~\ref{thmf}) 
generalizes and uniformizes the proofs of paper \cite{CPSZ}.

\begin{dfn}\label{celfib}
Let $X$ be a smooth projective variety over a field $F$.
We say a smooth projective morphism $f\colon Y\to X$ is a 
{\em cellular fibration} if it is a locally trivial fibration
whose fiber $\mathcal{F}$ is cellular, i.e.,
has a decomposition into affine cells (see \cite[\S 66]{EKM}). 
\end{dfn}

\begin{lem}\label{edigrah} 
Let $f\colon Y\to X$ be a cellular fibration.
Then $\M(Y)$ is isomorphic to 
$\M(X)\otimes \M(\mathcal{F})$. 
\end{lem}
\begin{proof}
We follow the proof of \cite[Proposition~1]{EG97}. 
Define the morphism
$$
\varphi\colon \bigoplus_{i\in\mathcal{I}} \M(X)(\codim B_i) \to \M(Y)
$$
to be the direct sum $\varphi=\bigoplus_{i\in\mathcal{I}}\varphi_i$,
where each $\varphi_i$ is given by the cycle 
$[\pr_Y^*(B_i)\cdot \Gamma_f]\in\CH(X\times Y)$ 
produced from  the graph cycle $\Gamma_f$ and 
the chosen (non-canonical) basis $\{B_i\}_{i \in \mathcal{I}}$
of $\CH(Y)$ over $\CH(X)$.
The realization of $\varphi$ coincides exactly with
the isomorphism of abelian groups $\CH(X)\otimes\CH(\mathcal{F})\to 
\CH(Y)$ constructed in \cite[Proposition~1]{EG97}.
Then by Manin's identity principle (see \cite[\S3]{Ma68}) 
$\varphi$ is an isomorphism.
\end{proof}

\begin{lem}\label{lemhomog}
Let $G$ be a linear algebraic group over a field $F$, 
$X$ be a projective homogeneous $G$-variety and $Y$ be a $G$-variety. 
Let $f\colon Y\to X$ be a $G$-equivariant projective morphism. 
Assume that the fiber of $f$ over $F(X)$ is 
isomorphic to $\mathcal{F}_{F(X)}$ for some variety $\mathcal{F}$ over $F$. 
Then $f$ is a locally trivial fibration with fiber $\mathcal{F}$.
\end{lem}

\begin{proof}
By the assumptions, we have $Y\times_X\Spec F(X)\simeq(\mathcal{F}\times 
X)\times_X\Spec F(X)$ as schemes over $F(X)$. Since $F(X)$ is a direct limit 
of $\mathcal{O}(U)$ taken over all non-empty affine open subsets $U$ of $X$, 
by [EGA IV, Corollaire~8.8.2.5] there exists $U$ such that $f^{-1}
(U)=Y\times_X U$ is isomorphic to $(\mathcal{F}\times X)\times_X 
U\simeq\mathcal{F}\times U$ as a scheme over $U$. Since $G$ acts transitively 
on $X$ and $f$ is $G$-equivariant, the map $f$ is a locally trivial 
fibration.
\end{proof}

\begin{cor}\label{corprod}
Let $X$ be a projective homogeneous $G$-variety, $Y$ be a projective variety 
such that $Y_{F(X)}\simeq\mathcal{F}_{F(X)}$ for some variety $\mathcal{F}$. 
Then the projection map $X\times Y\to X$ is a locally trivial fibration with 
fiber $\mathcal{F}$.
Moreover, if $\mathcal{F}$ is cellular, 
then $\M(X\times Y)\simeq\M(X)\otimes\M(\mathcal{F})$.
\end{cor}
\begin{proof}
Apply Lemma~\ref{lemhomog} to the projection map $X\times Y\to X$
and use Lemma~\ref{edigrah}.
\end{proof}

\begin{lem}\label{mainthmf}
Let $G$ be a semisimple linear algebraic group over $F$, $X$ and $Y$ be 
projective homogeneous $G$-varieties 
corresponding to parabolic subgroups $P$ 
and $Q$ of the split form $G_0$, $Q\subseteq P$. 
Denote by $f\colon Y\to X$ the map induced by the quotient map. 
If $G$ splits over $F(X)$ then $f$ is a cellular fibration with fiber
$\mathcal{F}=P/Q$.
\end{lem}

\begin{proof}
Since $G$ splits over $F(X)$, the fiber of $f$ over $F(X)$ is isomorphic to 
$(P/Q)_{F(X)}=\mathcal{F}_{F(X)}$. Now apply Lemma~\ref{lemhomog} and note 
that $\mathcal{F}$ is cellular.
\end{proof}

\begin{ex}\label{intmot} 
Let $P=P_\Theta$ be the standard parabolic subgroup of a split simple group 
$G_0$, corresponding to a subset $\Theta$ of the respective 
Dynkin diagram $\mathcal{D}$ 
(our enumeration of roots follows Bourbaki). In this 
notation the Borel subgroup corresponds to the empty set. Let $\xi$ be a 
cocycle in $H^1(F,G_0)$. 
Set $G={}_\xi G_0$ and $X={}_\xi(G_0/P)$. 
In particular, $G$ is a group of inner type
and $X$ is the respective projective homogeneous $G$-variety.
Denote by $q$ the degree of a splitting field of $G$ and 
by $d$ the index of the associated Tits
algebra (see \cite[Table~II]{Ti66}). For groups of type $\D_n$, 
we set $d$ to be the index 
of the Tits algebra associated with the vector representation.

Analyzing Tits indices of $G$
we see that {\em $G$ becomes 
split over $F(X)$} 
if the subset 
$\mathcal D\setminus\Theta$ contains one of the following vertices $k$ 
(cf. \cite[\S 7]{KR94}): 

\medskip

\noindent
\begin{tabular}{l||l|l|l|l|l}
$G_0$ & ${}^1\A_n$ & $\B_n$            & $\C_n$ & ${}^1\D_n$ & $\G$ \\ \hline
$k$ & $\gcd(k,d)=1$ & $k=n$;     & $k$ is odd; & $k=n-1,n$ if $d=1$; & any \\
&  & any $k$ in the &  & any $k$ in the &\\
& & Pfister case &  &  Pfister case &
\end{tabular}
 
\medskip

\noindent
\begin{tabular}{l||l|l|l|l}
$G_0$ & $\F$ & ${}^1\E_6$ & $\E_7$ & $\E_8$ \\\hline
$k$ & $k=1,2,3$; & $k=3,5$; & $k=2,5$; & $k=2,3,4,5$; \\
& any $k$ if  & $k=2,4$ if $d=1$; & $k=3,4$ if $d=1$; & any $k$ if \\
& $q=3$ & $k=1,6$ if $q$ is odd & $k\neq 7$ if $q=3$ & $q=5$\\
\end{tabular}

\medskip

\noindent
(here by the Pfister case we mean the case when the cocycle $\xi$ corresponds
to a Pfister form or its maximal neighbor)
\end{ex}

Case-by-case arguments of  paper \cite{CPSZ} show that 
under certain conditions the Chow motive
of a twisted flag variety $X$ can be expressed
in terms of the motive of a minimal flag.
These conditions cover almost all twisted flag varieties 
corresponding to groups of types $\A_n$ and $\B_n$ 
together with some examples of types $\C_n$, $\G$ and $\F$.
The following theorem together with Table~\ref{intmot}
provides a uniform proof of these results
and extends them to some other types. 

\begin{thm}\label{thmf} Let $Y$ and $X$ be taken as in Lemma~\ref{mainthmf}.
Then the Chow motive $\M(Y)$ of $Y$ is isomorphic
to a direct sum of twisted copies of the motive $\M(X)$, i.e.,
$$
\M(Y) \simeq \bigoplus_{i\ge 0} \M(X)(i)^{\oplus c_i},
$$
where $\sum c_it^i=P(\CH_*(\overline Y),t)/P(\CH_*(\overline X),t)$ is the
quotient of the respective Poincar\'e polynomials.
\end{thm}
\begin{proof} Apply Lemmas~\ref{mainthmf} and \ref{edigrah}.
\end{proof}

\begin{rem}\label{expldec} 
An explicit formula for $P(\CH_*(\overline X),t)$ involves degrees
of the basic polynomial invariants of $G_0$ and is provided in
\cite[Ch.~IV, Cor.~4.5]{Hi82}.
\end{rem}


\section{$J$-invariant and its properties}\label{JINV}

Fix a prime integer $p$.
To simplify the notation
we denote by $\Ch(X)$ the
Chow ring of a variety $X$ with $\zz/p$-coefficients and by 
$\ChO(X)$ the image of the restriction map 
$\CH(X;\zz/p)\to \CH(\overline X;\zz/p)$.

\begin{ntt}\label{sdef}
Let $G_0$ be a split semisimple linear algebraic group 
over a field $F$ with a 
split maximal torus $T$ and 
a Borel subgroup $B$ containing $T$.
Let $G={}_\xi G_0$ be a twisted form of $G_0$ 
given by a cocycle $\xi\in H^1(F,G_0)$.
Let $\XX={}_\xi(G_0/B)$ be the corresponding \emph{variety of complete flags}. 
Observe that the group $G$ splits over any field $K$ 
over which $\XX$ has a rational point, in particular, over
the function field $F(\XX)$. 
According to \cite{De74} the Chow ring 
$\Ch(\BXX)$ can be expressed in purely combinatorial terms and, 
therefore, depends only on the type of $G$ but not on the base field $F$.
\end{ntt}

\begin{ntt}\label{notgen} 
Let $\hat T$ denote the group of characters of $T$ and $S(\hat T)$
be the symmetric algebra. 
By $R$ we denote the image of the characteristic map 
$c\colon S(\hat T) \to \Ch(\BXX)$
(see \cite[(4.1)]{Gr58}).
According to \cite[Thm.6.4]{KM05} there is an embedding
\begin{equation}\label{ratemb}
R\subseteq \ChO(\XX),
\end{equation}
where the equality holds if the cocycle $\xi$ 
corresponds to a generic torsor.
\end{ntt}

\begin{ntt}
Let $\Ch(\BG)$ denote the Chow ring with $\zz/p$-coefficients
of the group $G_0$ over a splitting field of $\BXX$.
Consider the pull-back induced by the quotient map
$$
\pi\colon \Ch(\BXX)\to \Ch(\BG)
$$
According to \cite[p.~21, Rem.~$2^\circ$]{Gr58} 
$\pi$ is surjective with its kernel
generated by $R^+$, 
where $R^+$ stands for the subgroup of the non-constant elements of $R$.
\end{ntt}

\begin{ntt}
An explicit presentation of $\Ch(\BG)$  
is known for all types of $G$ and all torsion primes $p$ of $G$
(see \cite[Definition~3]{Gr58}). 
Namely, by \cite[Theorem~3]{Kc85} it is a quotient of the polynomial ring
in $r$ variables $x_1,\ldots,x_r$ of codimensions 
$d_1\le d_2\le\ldots\le d_r$ coprime to $p$, 
modulo an ideal generated by certain $p$-powers 
$x_1^{p^{k_1}},\ldots,x_r^{p^{k_r}}$ ($k_i\ge 0$, $i=1,\ldots,r$)
\begin{equation}\label{kcpres}
\Ch^*(\BG)=(\zz/p)[x_1,\ldots,x_r]/(x_1^{p^{k_1}},\ldots,x_r^{p^{k_r}}).
\end{equation}
In the case where $p$ is not a torsion prime of $G$ we have $\Ch^*(\BG)=\zz/p$, i.e.,
$r=0$.

Note that a complete list of
numbers $\{d_ip^{k_i}\}_{i=1\ldots r}$ 
called \emph{$p$-exceptional degrees} of $G_0$ 
is provided in \cite[Table II]{Kc85}.  
Taking the $p$-primary and $p$-coprimary parts of each $p$-exceptional degree 
from this table one restores the respective $k_i$ and $d_i$.
\end{ntt}

\begin{ntt} We introduce two orders on the set of additive generators 
of $\Ch(\BG)$, i.e., on the monomials 
$x_1^{m_1}\ldots x_r^{m_r}$. 
To simplify the notation, we will denote
the monomial $x_1^{m_1}\ldots x_r^{m_r}$ by $x^M$, where $M$ is an $r$-tuple
of integers $(m_1,\ldots,m_r)$. The codimension of $x^M$ will be denoted
by $|M|$. Observe that $|M|=\sum_{i=1}^r d_im_i$.
\begin{itemize}
\item
Given two $r$-tuples $M=(m_1,\ldots,m_r)$
and $N=(n_1,\ldots,n_r)$
we say $x^M\preccurlyeq x^N$ (or equivalently $M\preccurlyeq N$)
if $m_i\le n_i$ for all $i$. This gives a partial ordering on the
set of all monomials ($r$-tuples).
\item
Given two $r$-tuples $M=(m_1,\ldots,m_r)$
and $N=(n_1,\ldots,n_r)$ we say
$x^M\le x^N$ (or equivalently $M\le N$) if either $|M|<|N|$, or $|M|=|N|$ and 
$m_i\le n_i$ for the greatest $i$ such that $m_i\ne n_i$.
This gives a well-ordering on the set of all monomials ($r$-tuples)
known as the {\em DegLex} order. 
\end{itemize}
\end{ntt}

Now we are ready to give the main definition of the present paper.

\begin{dfn}\label{jinv} Let $G={}_\xi G_0$ 
be the twisted form of a split semisimple algebraic group $G_0$
over a field $F$ by means of a cocycle $\xi\in H^1(F,G_0)$
and $\XX={}_\xi (G_0/B)$ the respective variety of complete flags.
Let $\ChGi$ denote the image of the composite
$$
\Ch(\XX)\xra{\res} \Ch(\BXX) \xra{\pi} \Ch(\BG)
$$
Since both maps are ring homomorphisms, $\ChGi$ is a subring of 
$\Ch(\BG)$.

For each $1\le i\le r$ set $j_i$ to be the smallest non-negative
integer such that the subring $\ChGi$ contains an element $a$ 
with the greatest monomial $x_i^{p^{j_i}}$ 
with respect to the {\em DegLex} order on $\Ch(\BG)$, i.e.,
of the form 
$$
a=x_i^{p^{j_i}}+\sum_{x^M\lneq x_i^{p^{j_i}}} c_M x^M, \quad c_M\in\zz/p.
$$
The $r$-tuple of integers $(j_1,\ldots,j_r)$ will be called
the {\em $J$-invariant of $G$ modulo $p$} and will be denoted by $J_p(G)$.
\end{dfn}

Observe that if the Chow ring 
$\Ch(\BG)$ has only one generator, i.e., $r=1$, then
the $J$-invariant is equal to the smallest
non-negative integer $j_1$ such that $x_1^{p^{j_1}}\in \ChGi$.

\begin{ex} 
From the definition it follows that $J_p(G_E)\preccurlyeq J_p(G)$
for any field extension $E/F$.
Moreover, $J_p(G)\preccurlyeq (k_1,\ldots,k_r)$ by \eqref{kcpres}. 

According to \eqref{ratemb} the $J$-invariant
takes its maximal possible value $J_p(G)=(k_1,\ldots,k_r)$ if
the cocycle $\xi$ corresponds to a generic torsor.
Later on (see Corollary~\ref{trivJ}) it will
be shown that the $J$-invariant takes its minimal value 
$J_p(G)=(0,\ldots,0)$ if and only if 
the group $G$ splits by a finite field extension of degree coprime to $p$.
\end{ex}

The next example explains the terminology `$J$-invariant'.

\begin{ex}
Let $\phi$ be a quadratic form with trivial discriminant. 
In \cite[Definition~5.11]{Vi05} 
A.~Vishik introduced the notion of 
the $J$-invariant of $\phi$, a tuple of integers  
which describes the subgroup of rational
cycles on the respective maximal orthogonal Grassmannian.
This invariant provides an important tool for study of
algebraic cycles on quadrics. In particular, it was
one of the main ingredients used by A.~Vishik in his significant
progress on the solution of
Kaplansky's Problem. More precisely, 
in the notation of paper \cite{Vi06} the $J$-invariant of a quadric
corresponds to the upper row of its {\em elementary discrete invariant} 
(see \cite[Definition~2.2]{Vi06}). 

An equivalent but `dual' 
(in terms of non-rationality of cycles) definition of 
$J(\phi)$ was provided in \cite[\S~88]{EKM}. 
Using Theorem~\ref{thmf} one can show that $J(\phi)$ introduced in \cite{EKM}
can be expressed in terms of $J_2(\OO(\phi))=(j_1,\ldots,j_r)$ as follows:
$$
J(\phi)=\{2^ld_i\mid i=1,\ldots,r,\,0\le l\le j_i-1\}.
$$
Since all $d_i$ are odd, $J_2(\OO(\phi))$ 
is uniquely determined by $J(\phi)$.
\end{ex}

We have the following reduction formula 
(cf. \cite[Cor.~88.7]{EKM} in the case of quadrics).

\begin{prop}\label{thmKarp}
Let $G$ be a semisimple group of inner type
over a field $F$ and  
$\XX$ the variety of complete $G$-flags. 
Let $Y$ be a projective variety such that the map
$\CH^l(Y)\to\CH^l(Y_{F(x)})$ is surjective for all $x\in X$ and $l\le n$.
Then $j_i(G)=j_i(G_{F(Y)})$ for all $i$ such that
$d_ip^{j_i(G_{F(Y)})}\le n$.
\end{prop}
\begin{proof}
By \cite[Lemma~88.5]{EKM} the map
$\CH^l(X)\to\CH^l(X_{F(Y)})$ is surjective for all $l\le n$ and, therefore
$j_i(G)\le j_i(G_{F(Y)})$. The converse inequality is obvious.
\end{proof}
\begin{cor}
$J_p(G)=J_p(G_{F(t)})$.
\end{cor}
\begin{proof}
Take $Y=\mathbb{P}^1$ and apply Proposition~\ref{thmKarp}.
\end{proof}

\begin{ntt}To find restrictions on the possible values of $J_p(G)$
we use Steenrod $p$-th power operations introduced by P.~Brosnan. 
Recall (see \cite{Br03}) that if the characteristic 
of the base field $F$ is different from $p$ then 
one can construct \emph{Steenrod $p$-th power operations}
$$
S^l\colon\Ch^*(X)\to\Ch^{*+l(p-1)}(X),\qquad l\ge 0
$$
such that $S^0=\id$, the restriction $S^l\vert_{\Ch^l(X)}$ 
coincides with taking to the $p$-th power, $S^l\vert_{\Ch^i(X)}=0$ for $l>i$,
and the total operation 
$S^\bullet=\sum_{l\ge 0}S^l$ 
is a homomorphism of $\zz/p$-algebras 
compatible with pull-backs. In particular,
Steenrod operations preserve rationality of cycles.

In the case of projective homogeneous varieties over the field of complex numbers 
$S^l$ is compatible with its topological counterparts: 
the reduced power operation $\mathcal{P}^l$ if $p\ne 2$ and 
the Steenrod square $Sq^{2l}$ if $p=2$ 
(over complex numbers $\Ch^*(\BXX)$ can be viewed as a subring of 
the singular cohomology $\HH_{sing}^{2*}(\BXX,\zz/p)$).
Moreover, $\Ch^*(\BG)$ 
may be identified with the image of the pull-back map 
$\HH^{2*}_{sing}(\BXX,\zz/p)\to \HH^{2*}_{sing}(\BG,\zz/p)$. 
An explicit description of this image and 
formulae describing the action of $\mathcal{P}^l$ and $Sq^{2l}$
on $\HH^*_{sing}(\BG,\zz/p)$ can be found in \cite{MT91} for exceptional groups
and in \cite{EKM} for classical groups.

The action of the Steenrod operations on $\Ch(\BXX)$ and on $\Ch(\BG)$
can be described in purely combinatorial terms (see \cite{DZ07}) and, hence,
doesn't depend on the choice of a base field $F$.
\end{ntt}

The following lemma provides an important technical tool
for computing possible values of the $J$-invariant of $G$.

\begin{lem}\label{strestr}
Assume that in $\Ch^*(\BG)$ we have $S^l(x_i)=x_m^{p^s}$, and for any $i'<i$ $S^l(x_{i'})<x_m^{p^s}$ with respect to the DegLex order. Then $j_m\le j_i+s$.
\end{lem}

\begin{proof}
By definition there exists a cycle $\alpha\in\ChO(\XX)$ 
such that the leading term of $\pi(\alpha)$ is $x_i^{p^{j_i}}$. 
For the total operation we have
$$
S(x_i^{p^{j_i}})=S(x_i)^{p^{j_i}}=S^0(x_i)^{p^{j_i}}+S^1(x_i)^{p^{j_i}}+
\ldots+S^{d_i}(x_i)^{p^{j_i}}.
$$
In particular, $S^{lp^{j_i}}(x_i^{p^{j_i}})=S^l(x_i)^{p^{j_i}}$.
Applying $S^{lp^{j_i}}$ to $\alpha$
we obtain a rational cycle whose image 
under $\pi$ has the leading term $x_m^{p^{j_i+s}}$.
\end{proof}

\begin{ntt}\label{bigtab}
We summarize information about restrictions on the $J$-invariant
into the following table 
(numbers $r$, $d_i$ and $k_i$ are taken from \cite[Table~II]{Kc85}).
Recall that $r$ is the number of generators of $\Ch^*(\BG)$, $d_i$
are their codimensions and $k_i$ define the $p$-power relations. 

\medskip
\noindent
\begin{tabular}{|l|l|l|l|l|l|}
\hline
$G_0$ & $p$ & $r$ & $d_i$ & $k_i$ & $j_i$\\\hline
$\SL_n/\mu_m,\,m\mid n$ & $p\mid m$ & $1$ & $1$ & $p^{k_1}\parallel n$ &  \\
$\PGSp_n,\,2\mid n$ & $2$ & $1$ & $1$ & $2^{k_1}\parallel n$ &  \\
$\SO_n$ & $2$ & $[\frac{n+1}{4}]$ & $2i-1$ & $[\log_2\frac{n-1}{2i-1}]$ & $j_i\ge j_{i+l}$ if $2\nmid{i-1\choose l},$\\
& & & & & $j_i\le j_{2i-1}+1$\\
$\Spin_n$ & $2$ & $[\frac{n-3}{4}]$ & $2i+1$ & $[\log_2\frac{n-1}{2i+1}]$ & $j_i\ge j_{i+l}$ if $2\nmid{i\choose l},$\\
& & & & & $j_i\le j_{2i}+1$\\
$\PGO_{2n}$, $n>1$ & $2$ & $[\frac{n+2}{2}]$ & $1,\,i=1$ & $2^{k_1}\parallel n$ & $j_i\ge j_{i+l}$ if $2\nmid{i-2\choose l},$ \\
&&& $2i-3,\,i\ge 2$&$[\log_2\frac{2n-1}{2i-3}]$&  $j_i\le j_{2i-2}+1$\\
$\Spin_{2n}^{\pm},\,2\mid n$ & $2$ & $\frac{n}{2}$ & $1,\,i=1$ &$2^{k_1}\parallel n$&$j_i\ge j_{i+l}$ if $2\nmid{i-1\choose l}$\\
&&& $2i-1,\,i\ge 2$& $[\log_2\frac{2n-1}{2i-1}]$ &  $j_i\le j_{2i-1}+1$\\
$\G,\,\F,\,\E_6$ & $2$ & $1$ & $3$ & $1$ &\\
$\F,\,\E_6^{sc},\,\E_7$ & $3$ & $1$ & $4$ & $1$ &\\
$\E_6^{ad}$ & $3$ & $2$ & $1,\,4$ & $2,\,1$ &\\
$\E_7^{sc}$ & $2$ & $3$ & $3,\,5,\,9$ & $1,\,1,\,1$ & $j_1\ge j_2\ge j_3$\\
$\E_7^{ad}$ & $2$ & $4$ & $1,\,3,\,5,\,9$ & $1,\,1,\,1,\,1$ & $j_2\ge j_3\ge j_4$\\
$\E_8$ & $2$ & $4$ & $3,\,5,\,9,\,15$ & $3,\,2,\,1,\,1$ & $j_1\ge j_2\ge j_3,$\\
& & & & & $j_1\le j_2+1,\,j_2\le j_3+1$\\
$\E_8$ & $3$ & $2$ & $4,\,10$ & $1,\,1$ & $j_1\ge j_2$\\
$\E_8$ & $5$ & $1$ & $6$ & $1$ &\\
\hline
\end{tabular}
\medskip

The last column of the table follows from Lemma~\ref{strestr} and, hence,
requires $\mathrm{char}\,(F)\neq p$ restriction. All other columns 
are taken directly from  \cite[Table~II]{Kc85} and are independent on
the characteristic of the base field.
\end{ntt}

\section{Motivic decompositions}\label{MDEC}
In the present section we prove the main result of this paper 
(Theorem~\ref{thmMain}).
In the beginning of this section we describe
a basis of the subring of rational cycles of $\Ch(\BXX\times\BXX)$,
where $\XX$ is the variety of complete flags. The key
results here are Propositions~\ref{pairclm} and \ref{allratX}.
As a consequence, we obtain a motivic decomposition of $\XX$
(Theorem~\ref{thmKac}) in terms of certain motive $\RR_p(G)$.
Then using motivic decompositions of cellular fibrations obtained above
(Theorem~\ref{thmf}) we generalize Theorem~\ref{thmKac}
to arbitrary generically split projective homogeneous varieties.
At the end we discuss some properties of the motives $\RR_p(G)$.

\begin{ntt}
We use the notation of the previous section.
Let $G$ be a semisimple group of inner type over $F$ and
$\XX$ the respective variety of complete flags.
Let $R\subseteq \ChO(\XX)$ be the image of the characteristic map.
Consider the quotient map $\pi\colon \Ch(\BXX)\to\Ch(\BG)$.
Fix preimages $e_i$ of $x_i$ in $\Ch(\BXX)$. 
For an $r$-tuple $M=(m_1,\ldots,m_r)$
set $e^M=\prod_{i=1}^r e_i^{m_i}$. 
Set $N=(p^{k_1}-1,\ldots,p^{k_r}-1)$ and
$d=\dim \XX-|N|$.
\end{ntt}

\begin{lem}\label{Grothlem}
The Chow ring $\Ch(\BXX)$ is a free $R$-module
with a basis $\{e^M\}$, $M\preccurlyeq N$.
\end{lem}

\begin{proof}
Note that the subgroup $R^+$ of the non-constant elements
of $R$ is a nilpotent ideal in $R$.
Applying the Nakayama Lemma we obtain that
$\{e^M\}$ generates $\Ch(\BXX)$.
By \cite[(2)]{Kc85} $\Ch(\BXX)$ is a free $R$-module, hence,
for the Poincar\'e polynomials we have
$$
P(\Ch^*(\BXX),t)=P(\Ch^*(\BG),t)\cdot P(R^*,t).$$
Substituting $t=1$ we obtain that 
$$
\rk \Ch(\BXX)=\rk \Ch(\BG)\cdot \rk R.
$$
To finish the proof  observe that 
$\rk \Ch(\BG)$ 
coincides with the number of generators $\{e^M\}$.
\end{proof}

\begin{prop}\label{pairclm}
The pairing $R\times R \to \zz/p$ given by 
$(\alpha,\beta)\mapsto \deg(e^N\alpha\beta)$ is non-degenerated, i.e.,
for any non-zero element $\alpha\in R$ 
there exists $\beta$ such that $\deg(e^N\alpha\beta)\neq 0$. 
\end{prop}
\begin{proof}
Choose a homogeneous basis of $\Ch(\BXX)$.
Let $\alpha^\vee$ be the Poincar\'e dual of $\alpha$ with respect to this
basis. 
By Lemma~\ref{Grothlem} $\Ch(\BXX)$ is a free $R$-module with the basis 
$\{e^M\}$, hence,
expanding $\alpha^\vee$
we obtain
$$
\alpha^\vee=\sum_{M\preccurlyeq N} e^M\beta_M, \text{ where } \beta_M\in R.
$$
Note that if $M\neq N$ then $\codim \alpha\beta_M > d$, therefore,
$\alpha\beta_M=0$. So we can set $\beta=\beta_N$.
\end{proof}

From now on we fix a homogeneous $\zz/p$-basis $\{\alpha_i\}$ of $R$ and the
dual basis $\{\alpha_i^{\#}\}$ 
with respect to the pairing introduced in Proposition~\ref{pairclm}.

\begin{cor} \label{cbasis} 
For $|M|\le |N|$ we have
$$
\deg(e^M\alpha_i\alpha_j^{\#})=\begin{cases} 1, & M=N \text{ and } i=j; \\
0, & \text{ otherwise}.
\end{cases}
$$
\end{cor}
\begin{proof} 
If $M=N$, then it follows from the definition of the dual basis.
Assume $|M|<|N|$.  
If $\deg(e^M\alpha_i\alpha_j^{\#})\neq 0$, 
then $\codim(\alpha_i\alpha_j^{\#})>d$,
in contradiction with the fact that $\alpha_i\alpha_j^{\#}\in R$. 
Hence, we are reduced to the case $M\neq N$ and $|M|= |N|$.
Since $|M|=|N|$, $\codim(\alpha_i\alpha_j^{\#})=d$ and, hence, 
$R^+\alpha_i\alpha_j^{\#}=0$.
On the other hand 
there exists $i$ such that $m_i\ge p^{k_i}$ and
$e^{p^{k_i}}\in \Ch(\BXX)\cdot R^+$. Hence, $e^M\alpha_i\alpha_j^{\#}=0$.
\end{proof}

\begin{dfn}\label{grfil}
Given two pairs $(L,l)$ and $(M,m)$, 
where $L,M$ are $r$-tuples and $l,m$ are integers, 
we say $(L,l)\le (M,m)$ if either $L\lneq M$, or in the case $L=M$ 
we have $l\le m$.
We introduce a filtration on the ring $\Ch(\BXX)$ as follows: 
\begin{quote}
The $(M,m)$-th term $\Ch(\BXX)_{M,m}$ of the filtration
is the $\zz/p$-subspace spanned by the elements 
$e^I\alpha$ with $I\le M$, $\alpha\in R$ homogeneous, $\codim\alpha\le m$.
\end{quote}
Define the \emph{associated graded ring} as follows:
$$
A^{*,*}=\bigoplus_{(M,m)}A^{M,m},
\text{ where }
A^{M,m}=\Ch(\BXX)_{M,m}/\bigcup_{(L,l)\lneq(M,m)}\Ch(\BXX)_{L,l}.
$$
By Lemma~\ref{Grothlem} if $M\preccurlyeq N$ the graded component $A^{M,m}$ 
consists of the classes of elements 
$e^M\alpha$ with $\alpha\in R$ and $\codim\alpha= m$.
In particular, $\rk A^{M,m}=\rk R^m$.
Comparing the ranks we see that
$A^{M,m}$ is trivial when $M\not\preccurlyeq N$. 

Consider the subring $\ChO(\XX)$ of rational cycles 
with the induced filtration. 
The associated graded subring will be denoted by 
$A^{*,*}_{rat}$.
From the definition of the $J$-invariant it follows that
the elements $e_i^{p^{j_i}}$, $i=1,\ldots, r$, belong to $A^{*,*}_{rat}$.

Similarly, we introduce a filtration 
on the ring $\Ch(\BXX\times\BXX)$ as follows:
\begin{quote} 
The $(M,m)$-th term of the filtration is the $\zz/p$-subspace spanned 
by the elements $e^I\alpha\times e^L\beta$ with $I+L\le M$, 
$\alpha,\beta\in R$ homogeneous and $\codim\alpha+\codim\beta\le m$.
\end{quote} 
The associated graded ring will be denoted by $B^{*,*}$.
By definition 
$B^{*,*}$ is isomorphic to the tensor
product of graded rings $A^{*,*}\otimes_{\zz/p} A^{*,*}$. 
The graded subring associated to $\ChO(\XX\times \XX)$ will be denoted by 
$B^{*,*}_{rat}$.
\end{dfn}

\begin{ntt}\label{keyobserv}
The key observation is that due to Corollary~\ref{cbasis} we have
$$
\Ch(\BXX\times\BXX)_{M,m}\circ\Ch(\BXX\times\BXX)_{L,l}\subset
\Ch(\BXX\times\BXX)_{M+L-N,m+l-d}\text{ and }
$$
$$
(\Ch(\BXX\times\BXX)_{M,m})_\star(\Ch(\BXX)_{L,l})\subset
\Ch(\BXX)_{M+L-N,m+l-d}
$$
and, therefore, we have a correctly defined composition law
$$
\circ\colon B^{M,m}\times B^{L,l}\to B^{M+L-N,m+l-d}
$$
and the realization map (see \ref{chmot})
$$
\star\colon B^{M,m}\times A^{L,l}\to A^{M+L-N,m+l-d}
$$
In particular, $B^{N+*,d+*}$ 
can be viewed as a graded ring with respect to the composition
and $(\alpha\circ\beta)_\star=\alpha_\star\circ\beta_\star$.
Note also that both operations preserve rationality of cycles.
\end{ntt}

The proof of  the following result is based on the fact that the group $G$
splits over $F(\XX)$.

\begin{lem}\label{sigrat}
The classes of the elements $e_i\times 1-1\times e_i$ in 
$B^{*,*}$, $i=1,\ldots,r$,  belong to $B^{*,*}_{rat}$.
\end{lem}

\begin{proof}  
 Fix an $i$.
Since $G$ splits over $F(\XX)$, $F(\XX)$ is a splitting field of $\XX$ and 
by Lemma~\ref{gensp}
there exists a cycle in $\overline{\Ch^{d_i}}(\XX\times \XX)$ of the form
$$
\xi=e_i\times 1+\sum_s \mu_s\times \nu_s+1\times\mu,
$$
where $\codim\mu_s,\codim\nu_s<d_i$. Then the cycle
$$
\pr_{13}^*(\xi)-\pr_{23}^*(\xi)=(e_i\times 1-1\times e_i)\times 1+
\sum_s (\mu_s\times 1-1\times\mu_s)\times\nu_s
$$
belongs to $\ChO(\XX\times \XX\times \XX)$, 
where $\pr_{ij}$ denotes the projection
on the product of the $i$-th and $j$-th factors. 
Applying Corollary~\ref{corprod} 
to the projection 
$\pr_{12}\colon \XX\times \XX\times \XX\to \XX\times \XX$ 
we conclude that  
there exists a (non-canonical) $\Ch(\XX\times \XX)$-linear
isomorphism\ $\Ch(\XX\times \XX\times \XX)\simeq
\Ch(\XX\times \XX)\otimes\Ch(\BXX)$, where
$\Ch(\XX\times \XX)$ acts on the left-hand side via $\pr_{12}^*$. 
This gives rise to
a $\Ch(\XX\times \XX)$-linear retraction $\delta$ to the pull-back map
$\pr_{12}^*$.
Since the construction of the retraction preserves base change,
it preserves rationality of cycles.
Hence, passing to a splitting field we obtain a rational cycle
$$
\bar\delta(\pr_{13}^*(\xi)-\pr_{23}^*(\xi))=e_i\times 1-1\times e_i+
\sum_s (\mu_s\times 1-1\times\mu_s)\bar\delta(1\times1\times\nu_s)
$$
whose image in $B^{*,*}_{rat}$ is $e_i\times 1-1\times e_i$.
\end{proof}

We will write $(e\times 1-1\times e)^M$ for the product
$\prod_{i=1}^r (e_i\times 1-1\times e_i)^{m_i}$
and ${M\choose L}$ for the product of binomial coefficients
$\prod_{i=1}^r{m_i\choose l_i}$. We assume that ${m_i\choose l_i}=0$ 
if $l_i>m_i$.
In the computations we will extensively use the following two
formulae
(the first follows directly from Corollary~\ref{cbasis} and the second one
is a well-known binomial identity).

\begin{ntt}\label{karpform}
Let $\alpha$ be an element of $R^*$ and $\alpha^\#$ be its dual
with respect to the non-degenerate pairing from~\ref{pairclm}, 
i.e. $\deg(e^N\alpha\alpha^\#)=1$. Then we have
$$
((e\times 1-1\times e)^M(\alpha^\#\times 1))_\star(e^L\alpha)=
{M \choose M+L-N}(-1)^{M+L-N}e^{M+L-N}.
$$
Indeed, expanding the brackets in the left-hand side,
we obtain
$$
\Big(\sum_{I\preccurlyeq M}(-1)^I{M\choose I}e^{M-I}\alpha^\#\times e^I\Big)_\star(e^L\alpha),
$$
and it remains to apply Corollary~\ref{cbasis}.
\end{ntt}

\begin{ntt}[Lucas' Theorem]\label{lucas}
The following identity holds
$$
{n\choose m}\equiv \prod_{i\ge 0} {n_i \choose m_i} \mod p,
$$
where $m=\sum_{i\ge 0}m_ip^i$ and $n=\sum_{i\ge 0} n_i p^i$ are
the base $p$ presentations of $m$ and $n$.
\end{ntt}

Let $J=J_p(G)=(j_1,\ldots,j_r)$ be the $J$-invariant of $G$
(see Definition~\ref{jinv}). Set $K=(k_1,\ldots,k_r)$.

\begin{prop}\label{allratX}
Let $\{\alpha_i\}$ be a homogeneous $\zz/p$-basis of $R$.
Then the set of elements $\mathcal{B}=\{e^{p^JL}\alpha_i\mid 
L\preccurlyeq p^{K-J}-1\}$ 
forms a $\zz/p$-basis of $A^{*,*}_{rat}$.
\end{prop}
\begin{proof}
According to Lemma~\ref{Grothlem} 
the elements from $\mathcal{B}$ are linearly independent. 
Assume $\mathcal{B}$ does not generate $A^{*,*}_{rat}$.
Choose an element $\omega\in A^{M,m}_{rat}$
of the smallest index $(M,m)$
which is not in the linear span of $\mathcal{B}$.
By definition of $A^{M,m}$ (see Definition~\ref{grfil}) 
$\omega$ can be written as 
$\omega=e^M\alpha$, where $M\preccurlyeq N$, $\alpha\in R^m$ and
$M$ can not be presented as $M=p^JL'$ for an $r$-tuple $L'$.
The latter means that in the decomposition of $M$ into $p$-primary
and $p$-coprimary components $M=p^SL$, where 
$M=(m_1,\ldots,m_r)$, $S=(s_1,\ldots,s_r)$, 
$L=(l_1,\ldots,l_r)$ and $p\nmid l_k$ 
for $k=1,\ldots,r$, 
we have 
$J\not\preccurlyeq S$.
Choose an $i$ such that $s_i<j_i$. Denote $M_i=(0,\ldots,0,m_i,0,\ldots,0)$
and $S_i=(0,\ldots,0,s_i,0,\ldots,0)$,
where $m_i$ and $s_i$ stand at the $i$-th place.

Set $T=N-M+M_i$.                    
By Lemma~\ref{sigrat} and \ref{karpform} together with
 observation~\ref{keyobserv}
the element
$$
((e\times 1-1\times e)^T(\alpha^\#\times 1))_\star(e^M\alpha)=
{p^{k_i}-1 \choose m_i}(-1)^{m_i}e^{m_i}
$$
belongs to $A^{M_i,0}_{rat}$. By \ref{lucas} we have
$p\nmid {{p^{k_i}-1}\choose m_i}$ and, therefore, this element is non-trivial.
Moreover, since $s_i<j_i$, this element is not in the span of $\mathcal{B}$.
Since $(M,m)$ was chosen to be the smallest index
and $(M_i,0)\le (M,m)$ we obtain that $(M,m)=(M_i,0)$.
Repeating the same arguments 
for $T=N-M_i+p^{S_i}$ we obtain that
$M_i=p^{S_i}$, i.e., $l_i=1$. 

Now let
$\gamma$ be a representative of $\omega=e_i^{p^{s_i}}$ in $\ChO(\XX)$. 
Then its image $\pi(\gamma)$ in $\ChGi$
has the leading term $x_i^{p^{s_i}}$ 
with $s_i<j_i$. This contradicts the definition of the $J$-invariant.
\end{proof}

\begin{cor}\label{allratXX}
The elements 
$$
\{(e\times 1-1\times e)^S(e^{p^JL}\alpha_i\times e^{p^J(p^{K-J}-1-M)}\alpha_j^{\#})\mid 
L, M\preccurlyeq p^{K-J}-1,\; S\preccurlyeq p^J-1\}
$$
form a $\zz/p$-basis of $B^{*,*}_{rat}$. In particular,
the ones such that  $S=p^J-1$ and $L=M$
form a basis of $B^{N,d}_{rat}$.
\end{cor}
\begin{proof}
According to Lemma~\ref{Grothlem} these elements are linearly independent 
and their number is $p^{|2K-J|}(\rk R)^2$. They are rational 
by Definition~\ref{grfil} and
Lemma~\ref{sigrat}.
Applying Corollary~\ref{corprod} to the projection $\BXX\times \BXX\to\BXX$
we obtain that
$$
\rk B^{*,*}_{rat} =
\rk \ChO(\XX\times \XX)= 
\rk\ChO(\XX)\cdot\rk\Ch(\BXX),
$$
where the latter coincides with $\rk A^{*,*}_{rat}\cdot p^{|K|}\rk R =
 p^{|2K-J|}(\rk R)^2$ by Lemma~\ref{Grothlem} and Proposition~\ref{allratX}.
\end{proof}

\begin{lem}\label{clidem}
The elements
$$
\theta_{L,M,i,j}=
(e\times 1-1\times e)^{p^J-1}
(e^{p^JL}\alpha_i\times e^{p^J(p^{K-J}-1-M)}\alpha_j^{\#}),\; 
L,M\preccurlyeq p^{K-J}-1,
$$
belong to $B^{*,*}_{rat}$ and
satisfy the relations
$\theta_{L,M,i,j}\circ\theta_{L',M',i',j'}=\delta_{LM'}\delta_{ij'}
\theta_{L',M,i',j}$  and
$\sum_{L,i}\theta_{L,L,i,i}=\Delta_{\BXX}$.
\end{lem}

\begin{proof}
Expanding the brackets and using the identity ${p^j-1\choose i}\equiv (-1)^i\mod p$,
we see that
$$
\theta_{L,M,i,j}=\sum_{I\preccurlyeq p^J-1}e^{p^JL+I}\alpha_i\times e^{N-p^JM-I}\alpha_j^\#,
$$
and the composition relation follows from Corollary~\ref{cbasis}.
By definition we have
$$
\sum_{L,i}\theta_{L,L,i,i}=\sum_{I\preccurlyeq N,i}e^I\alpha_i\times e^{N-I}\alpha_i^\#.
$$
By Corollary~\ref{cbasis} 
the latter sum acts trivially on all basis elements of $\Ch(\BXX)$
and, hence, coincides with the diagonal.
\end{proof}

We are now ready to provide a motivic decomposition of the variety
of complete flags.

\begin{thm}\label{thmKac}
Let $G$ be a semisimple linear algebraic group of inner type over a field $F$
and $\XX$ be the variety of complete $G$-flags.
Let $p$ be a prime.
Assume that $J_p(G)=(j_1,\ldots,j_r)$.
Then the motive of $\XX$ is isomorphic to the direct sum
$$
\M(\XX;\zz/p)\simeq \bigoplus_{i\ge 0}\RR_p(G)(i)^{\oplus c_i}, 
$$
where the motive $\RR_p(G)$ is indecomposable, 
its Poincar\'e polynomial over a splitting field is given by
\begin{equation}\label{poincr}
P(\overline{\RR_p(G)},t)=\prod_{i=1}^r \frac{1-t^{d_ip^{j_i}}}{1-t^{d_i}},
\end{equation}
and the integers $c_i$ are the coefficients of the quotient 
$$
\sum_{i\ge 0} c_i t^i=P(\Ch^*(\BXX),t)/P(\overline{\RR_p(G)},t).
$$
\end{thm}

\begin{proof}
Consider the projection map
$$
f^0\colon\ChO(\XX\times \XX)_{N,d}\to B^{N,d}_{rat}.
$$
Observe that the kernel of $f^0$ is nilpotent. Indeed, any element
$\xi$ from $\ker f^0$ belongs to $\Ch(\BXX\times\BXX)_{M,m}$ for
some $(M,m)\lneq (N,d)$ which depends on $\xi$. 
Then by~\ref{keyobserv} its $i$-th composition power
$\xi^{\circ i}$ belongs to the graded component
$\Ch(\BXX\times\BXX)_{iM-(i-1)N,im-(i-1)d}$,
and, therefore, becomes trivial for $i$ big enough.

By Lemma~\ref{clidem}  the elements $\theta_{L,L,i,j}$ 
form a family of pairwise-orthogonal
idempotents whose sum is the identity. Therefore, by Proposition~\ref{nilpcor} 
there exist pair-wise orthogonal idempotents $\varphi_{L,i}$ in 
$\ChO(\XX\times \XX)$ which are mapped to $\theta_{L,L,i,i}$ 
and whose sum is the identity. 

Recall (see \ref{chmot}) that
given two correspondences $\phi$ and $\psi$ in $\ChO(\XX\times \XX)$
of degrees $c$ and $c'$ respectively 
its composite $\phi\circ\psi$ has degree $c+c'$.
Using this fact we conclude that 
the homogeneous components of $\varphi_{L,i}$ of codimension $\dim \XX$ 
are pair-wise
orthogonal idempotents whose sum is the identity. Hence, we may
assume that 
$\varphi_{L,i}$ belong to $\overline{\Ch^{\dim \XX}}(\XX\times \XX)$.

We now show that $\varphi_{L,i}$ are indecomposable. 
By Corollary~\ref{allratXX} and Lemma~\ref{clidem} 
the ring $(B^{N,d}_{rat},\circ)$ can be identified
with a product of matrix rings over $\zz/p$
$$
B^{N,d}_{rat}\simeq\prod_{s=0}^d\End((\zz/p)^{{p^{|K-J|}\rk R^s}}).
$$
By means of this identification $\theta_{L,L,i,i}\colon 
e^{p^JM}\alpha_j\mapsto \delta_{L,M}\delta_{i,j}e^{p^JL}\alpha_i$ 
is an idempotent of rank $1$ and, therefore, is indecomposable. 
Since the kernel of $f^0$ is nilpotent, 
the $\varphi_{L,i}$ are indecomposable as well.

Next we show that $\varphi_{L,i}$ is isomorphic to $\varphi_{M,j}$. 
In the ring 
$B^{*,*}_{rat}$ mutually inverse isomorphisms between them are given by
$\theta_{L,M,i,j}$ and $\theta_{M,L,j,i}$. Let
$$
f\colon\ChO(\XX\times \XX)\to B^{*,*}_{rat}
$$
be the \emph{leading term} map; it means that for any $\gamma\in
\ChO(\XX\times \XX)$ we find the smallest degree $(I,s)$ such that
$\gamma$ belongs to $\ChO(\XX\times \XX)_{I,s}$ and set $f(\gamma)$ to be
the image of $\gamma$ in $B^{I,s}_{rat}$. 
Note that $f$ is not a homomorphism 
but satisfies the condition that $f(\xi)\circ f(\eta)$ equals either 
$f(\xi\circ\eta)$ or $0$. Choose preimages $\psi_{L,M,i,j}$ and
$\psi_{M,L,j,i}$ of $\theta_{L,M,i,j}$ and $\theta_{M,L,j,i}$ 
by means of $f$. 
Applying Lemma~\ref{izvrat} we obtain mutually inverse isomorphisms 
$\vartheta_{L,M,i,j}$ and $\vartheta_{M,L,j,i}$ between 
$\varphi_{L,i}$ and 
$\varphi_{M,j}$. By the definition of $f$ it remains to take their homogeneous components of the 
appropriate degrees.

Applying now Lemma~\ref{rostnil} and Corollary~\ref{rostcor} 
to the restriction map
$$
\res_{F}\colon\End(\M(\XX;\zz/p))\to\overline{\End}(\M(\XX;\zz/p))
$$
and the family of idempotents $\varphi_{L,i}$
we obtain a family of pair-wise orthogonal idempotents
$\phi_{L,i}\in\End(\M(\XX;\zz/p))$ 
such that 
$$
\Delta_{\XX}=\sum_{L,i} \phi_{L,i}.
$$
Since $\res_{F_s/F}$ 
lifts isomorphisms, for the respective motives we have
$(\XX,\phi_{L,i})\simeq (\XX,\phi_{0,0})(|L|+\codim \alpha_i)$ 
for all $L$ and $i$ 
(see \ref{motex}). 
The twists $|L|+\codim \alpha_i$ can be easily recovered from
the explicit formula for $\theta_{L,L,i,i}$ (see Lemma~\ref{clidem}).
Denoting $\RR_p(G)=(\XX,\phi_{0,0})$
we obtain the desired motivic decomposition.

Finally, consider the motive $\RR_p(G)$ over a splitting field.
The idempotent $\theta_{0,0,0,0}$ splits into the sum of 
pair-wise orthogonal (non-rational) 
idempotents $e^I\times e^{N-I}1^\#$, $I\preccurlyeq p^J-1$.
The motive corresponding to each summand is 
isomorphic to $(\zz/p)(|I|)$. Therefore, we obtain the decomposition
into Tate motives
$$
\overline{\RR_p(G)}\simeq\bigoplus_{I\preccurlyeq p^J-1}(\zz/p)(|I|),
$$
which gives formula~\eqref{poincr} for the Poincar\'e polynomial.
\end{proof}

As a direct consequence of the proof we obtain

\begin{cor}\label{remKrull}
Any direct summand of $\M(\XX;\zz/p)$ 
is isomorphic to a direct sum of twisted copies of $\RR_p(G)$.
\end{cor}
\begin{proof}
Indeed, in the ring $B^{N,d}_{rat}$ any idempotent is isomorphic 
to a sum of idempotents $\theta_{L,L,i,i}$, 
and the map $f^0$ lifts isomorphisms. 
\end{proof}

\begin{rem}
Corollary~\ref{remKrull} can be viewed as a particular
case of the Krull-Schmidt Theorem proven by V.~Chernousov
and A.~Merkurjev (see \cite[Corollary~9.7]{CM06}).
\end{rem}

\begin{dfn} Let $G$ be a linear algebraic group over a field $F$
and $X$ a projective homogeneous $G$-variety. We say $X$ is 
{\em generically split} if the group $G$ splits over the generic
point of $X$.
\end{dfn}

The main result of the present paper is the following

\begin{thm}\label{thmMain}
Let $G$ be a semisimple linear algebraic group of inner type
over a field $F$ and $p$ be a prime integer.
Let $X$ be a generically split
projective homogeneous $G$-variety. 
Then the motive of $X$ with $\zz/p$-coefficients is 
isomorphic to the direct sum
$$
\M(X;\zz/p)\simeq \bigoplus_{i\ge 0}\RR_p(G)(i)^{\oplus a_i},
$$
where $\RR_p(G)$ is an indecomposable motive; 
Poincar\'e polynomial $P(\overline{\RR_p(G)},t)$ 
is given by \eqref{poincr} and,
hence, only depends on the $J$-invariant of $G$;
the $a_i$'s are the coefficients of the quotient polynomial 
$$
\sum_{i\ge 0} a_i t^i=P(\CH^*(\BX),t)/P(\overline{\RR_p(G)},t).
$$
\end{thm}
\begin{proof}
Let $\XX$ be the variety of complete $G$-flags.
According to Theorem~\ref{thmf} the motive of $Y=\XX$ is isomorphic to
a direct sum of twisted copies of the motive of $X$.
To finish the proof we apply Theorem~\ref{thmKac} and 
Corollary~\ref{remKrull}.
\end{proof}

We now provide several
properties of
$\RR_p(G)$ which will be extensively used in the applications.

\begin{prop}\label{thmComp}
Let $G$ and $G'$ be two semisimple algebraic 
groups of inner type over $F$, $\XX$ and $\XX'$ be 
the corresponding varieties of complete flags.
\begin{itemize}
\item[(i)] {\rm\bf (base change)} For any field extension $E/F$ we have
$$\RR_p(G)_E\simeq\bigoplus_{i\ge 0}\RR_p(G_E)(i)^{\oplus a_i},$$
where $\sum a_it^i=P(\overline{\RR_p(G)},t)/P(\overline{\RR_p(G_E)},t)$.
\item[(ii)] {\rm\bf (transfer argument)} If $E/F$ is a field extension of degree 
coprime to $p$ then $J_p(G_E)=J_p(G)$ and $\RR_p(G_E)=\RR_p(G)_E$. Moreover, 
if $\RR_p(G_E)\simeq\RR_p(G'_E)$ then $\RR_p(G)\simeq\RR_p(G')$.
\item[(iii)] {\rm\bf (comparison lemma)} If $G$ splits over $F(\XX')$ 
and $G'$ splits 
over $F(\XX)$ then $\RR_p(G)\simeq\RR_p(G')$.
\end{itemize}
\end{prop}
\begin{proof}
The first claim follows from Theorem~\ref{thmKac} and 
Corollary~\ref{remKrull}. 
To prove the second claim note that 
$E$ is rank preserving with respect to 
$\XX$ and $\XX\times \XX$ by Lemma~\ref{surj}. 
Now $J_p(G_E)=J_p(G)$ by 
Lemma~\ref{prat}, and hence $\RR_p(G_E)=\RR_p(G)_E$ 
by the first claim. The 
remaining part of the claim follows 
from Corollary~\ref{transf} applied to the 
variety $\XX\coprod \XX'$.

We now prove the last claim.
The variety $\XX\times \XX'$ is the variety of 
complete $G\times G'$-flags. By Corollary~\ref{corprod} applied to the
projections $\XX\times \XX'\to \XX$ and 
$\XX\times \XX'\to \XX'$ we can express 
$\M(\XX\times \XX';\zz/p)$ in terms of $\M(\XX;\zz/p)$ and $\M(\XX';\zz/p)$.
The latter motives can be expressed in terms of $\RR_p(G)$ and $\RR_p(G')$.
Now the claim follows from the Krull-Schmidt theorem 
(see Corollary~\ref{remKrull}).
\end{proof}

\begin{cor}
We have $\RR_p(G)\simeq\RR_p(G_{an})$, where $G_{an}$ is the semisimple
anisotropic kernel of $G$.
\end{cor}

Finally, we provide conditions
which allow to lift a motivic decomposition
of a generically split homogeneous variety with $\zz/m$-coefficients
to a decomposition with $\zz$-coefficients.

\begin{ntt}
Let $m$ be a positive integer.
We say a polynomial $g(t)$ is $m$-{\it positive}, 
if $g\neq 0$, $P(\overline{\RR_p(G)},t)\mid g(t)$ 
and the quotient polynomial $g(t)/P(\overline{\RR_p(G)},t)$
has non-negative coefficients
for all primes $p$ dividing $m$.
\end{ntt}

\begin{prop}\label{intd}
Let $G$ be a semisimple linear algebraic group of inner type
over a field $F$ and $X$ be a generically split
projective homogeneous $G$-variety. Assume that 
$X$ splits by a field extension of degree $m$.
Let $f(t)$ be an $m$-positive polynomial dividing $P(\M(\BX),t)$ which
can not be presented as a sum of two $m$-positive polynomials.
Then the motive of $X$ with integer coefficients
splits as a direct sum
$$
\M(X;\zz)\simeq \bigoplus_{i}\RR_i(c_i),\quad c_i\in\zz,
$$
where $\RR_i$ are indecomposable and
$P(\overline{\RR}_i,t)=f(t)$ for all $i$.
Moreover, if $m=2,3,4$ or $6$, then
all motives $\RR_i$ are isomorphic up to twists.
\end{prop}
\begin{proof} First, we apply Corollary~\ref{karpcor}
to obtain a decomposition with $\zz/m$-coefficients. 
By Lemma~\ref{surj} our field extension is rank preserving so
we can apply Theorem~\ref{lem47}
to lift the decomposition to the category of motives with $\zz$-coefficients.
\end{proof}


\section{Applications of the $J_p$ and of the motive $\RR_p$}
\label{APPL}

Let $G$ be a semisimple 
group of inner type over $F$ and $\XX$ the variety of complete
$G$-flags.

First, we obtain
the following expression for the {\em canonical $p$-dimension}
of $\XX$ (see \cite[\S90]{EKM}).

\begin{prop} In the notation of Theorem~\ref{thmKac} we have
$$
\mathrm{cd}_p(\XX)=\sum_{i=1}^r d_i(p^{j_i}-1).
$$
\end{prop}
\begin{proof}
Follows from Proposition~\ref{allratX} and \cite[Theorem~5.8]{KM05}.
\end{proof}

Let $X$ be a smooth projective variety which has a splitting field.

\begin{lem}\label{trlem}
For any $\phi,\psi\in\CH^*(\BX\times\BX)$ one has
$$
\deg((\pr_2)_*(\phi\cdot\psi^t))=\tr((\phi\circ\psi)_\star).
$$
\end{lem}
\begin{proof}
Choose a homogeneous basis $\{e_i\}$ of $\CH^*(\BX)$. 
Let $\{e_i^\vee\}$ be its Poincar\'e dual. 
Since both sides of the relation under proof are bilinear, 
it suffices to check the assertion for 
$\phi=e_i\times e_j^\vee$ and $\psi=e_k\times e_l^\vee$. 
In this case both sides of the relation
are equal to $\delta_{il}\delta_{jk}$.
\end{proof}

Denote by $d(X)$ the greatest common divisor 
of the degrees of all zero cycles on $X$ and 
by $d_p(X)$ its $p$-primary component.

\begin{cor}\label{trcor} Let $m$ be an integer.
For any $\phi\in\CHO(X\times X;\zz/m)$ we have 
$$\gcd(d(X),m)\mid\tr(\phi_\star).$$
\end{cor}
\begin{proof}
Set $\psi=\Delta_{\BX}$ and apply Lemma~\ref{trlem}.
\end{proof}

\begin{cor}\label{degr}
Assume that $\M(X;\zz/p)$ has a direct summand $M$. Then
\begin{enumerate}
\item $d_p(X)\mid P(\overline{M},1)$;
\item if $d_p(X)=P(\overline{M},1)$ and the kernel of the restriction 
$\End(\M(X))\to\End(\M(\BX))$ consists of nilpotents, 
then $M$ is indecomposable.
\end{enumerate}
\end{cor}
\begin{proof}
Set $q=d_p(X)$ for brevity. Let
$M=(X,\phi)$. 
By Corollary~\ref{karpcor} there exists an idempotent
$\varphi\in\End(\M(X);\zz/q)$ such that $\varphi\mod p=\phi$. 
Then $\res(\varphi)\in\End(\M(\BX);\zz/q)$ is a rational idempotent. 
Since every projective module over $\zz/q$ is free, we have
$$
\tr(\res(\varphi)_\star)=\rk_{\zz/q}(\res(\varphi)_\star)=
\rk_{\zz/p}(\res(\phi)_\star)=P(\overline{M},1)\mod q,
$$
and the first claim follows from Corollary~\ref{trcor}. 
The second claim follows from the first one, 
since the second assumption implies 
that for any non-trivial direct summand $M'$ of $M$ 
we have $P(\overline{M'},1)<P(\overline{M},1)$.
\end{proof}

\begin{ntt}
Denote by $n(G)$ the greatest common divisor of degrees of all
finite splitting fields of $G$ and by $n_p(G)$ its $p$-primary component.
Note that $n(G)=d(\XX)$ and $n_p(G)=d_p(\XX)$.
\end{ntt}

We obtain the following estimate on $n_p(G)$ in terms of the $J$-invariant 
(cf. \cite[Prop.~88.11]{EKM} in the case of quadrics).
\begin{prop}\label{npG}
Let $G$ be a semisimple linear algebraic group of inner type with 
$J_p(G)=(j_1,\ldots,j_r)$. Then
$$
n_p(G)\le p^{\sum_i j_i}.
$$
\end{prop}
\begin{proof}
Follows from Theorem~\ref{thmKac} and Corollary~\ref{degr}.
\end{proof}

\begin{cor}\label{trivJ} The following statements are equivalent:
\begin{itemize}
\item $J_p(G)=(0,\ldots,0)$;
\item $n_p(G)=1$;
\item $\RR_p(G)=\zz/p$.
\end{itemize}
\end{cor}
\begin{proof}
If $J_p(G)=(0,\ldots,0)$ then $n_p(G)=1$ by 
Proposition~\ref{npG}. If $n_p(G)=1$ 
then there exists a splitting field $L$ of degree $m$ prime to $p$ and, 
therefore, $\RR_p(G)=\zz/p$ by transfer argument \ref{thmComp}(ii). 
The remaining implication is obvious.
\end{proof}


\section{Examples}\label{mainexam}

In the present section we provide
examples of motivic decompositions of projective homogeneous varieties
obtained using Theorem~\ref{thmMain}. 

\paragraph{The case $r=d_1=1$.}
According to Table~\ref{bigtab} this corresponds to the case when
$G$ is of type $\A_n$ or $\C_n$. Let
$A$ be a central simple algebra corresponding to $G$. 
We have $A=\mathrm{M}_m(D)$, 
where $D$ is a division algebra of index $d\ge 1$ over a field $F$. 
Let $p$ be a prime divisor of $d$. 
Observe that according to Table~\ref{bigtab} $J_p(G)=(j_1)$ for some
$j_1\ge 0$.
Let $X_\Theta$ be the projective homogeneous $G$-variety 
given by a subset $\Theta$ of vertices of the respective Dynkin diagram 
such that $p\nmid j$ for some $j\notin\Theta$ (cf.~Example~\ref{intmot}).
Then by Theorem~\ref{thmMain} we obtain that
\begin{equation}\label{sbdec}
\M(X_\Theta;\zz/p)\simeq\bigoplus_{i\ge 0}\RR_p(G)(i)^{\oplus a_i},
\end{equation}
where $\RR_p(G)$ is indecomposable and 
$$
\overline{\RR_p(G)}\simeq\bigoplus_{i=0}^{p^{j_1}-1}(\zz/p)(i).
$$

We now identify $\RR_p(G)$.
Using the comparison lemma (see Proposition~\ref{thmComp})
we conclude that 
$\RR_p(G)$ only depends on $D$, so we may assume $m=1$. 
By Table~\ref{bigtab} we have $p^{j_1}\mid d$, but on the other hand by Proposition~\ref{npG} we have
$n_p(G)\le p^{j_1}$. 
Therefore, $p^{j_1}$ is a $p$-primary part of $d$.
 
We have $D\simeq D_p\otimes_F D'$, where $p^{j_1}=\ind(D_p)$
and $p\nmid\ind(D')$. Passing to a splitting field of $D'$
of degree prime to $p$ and using Proposition~\ref{thmComp} 
we conclude that the motives of $X_\Theta$ 
and $\SB(D_p)$ are direct sums of twisted $\RR_p(G)$.
Comparing the Poincar\'e polynomials we conclude that 

\begin{lem}\label{idsb}
$\M(\SB(D_p);\zz/p)\simeq \RR_p(G)$.
\end{lem}

Applying Proposition~\ref{intd} to $X=\SB(D)$ and comparing
the Poincar\'e polynomials of $\M(X)$ and $\RR_i$ we obtain
that 
\begin{cor}
The motive of $\SB(D)$ with
integer coefficients is indecomposable.
\end{cor}

\begin{rem} Indeed, we provided
a uniform proof of the results
of paper \cite{Ka96}. Namely, the decomposition of $\M(\SB(A);\zz/p)$
(see \cite[Cor.~1.3.2]{Ka96})
and indecomposability of $\M(\SB(D);\zz)$ (see \cite[Thm.~2.2.1]{Ka96}).
\end{rem}

\paragraph{The case $r=1$ and $d_1>1$.} 
According
to Table~\ref{bigtab} this holds if
\begin{itemize}
\item[$p=2$:] $\G$, $\F$, $\E_6$ or 
$G$ is a strongly inner form of type $\B_3$, $\B_4$, $\D_4$, $\D_5$;
\item[$p=3$:] $G$ is a group of type $\F$, $\E_7$ or strongly
inner form of type $\E_6$;
\item[$p=5$:] $G$ is a group of type $\E_8$.
\end{itemize}
We say a group $G$ is {\em strongly inner} over a field $F$ if it is
the twisted form by means of a cocycle from $H^1(F,G_0)$, where
$G_0$ is the simply-connected split group over $F$ of the same type as $G$
(see \ref{sdef}). 

Observe that in all these cases $J_p(G)=(0)$ or $(1)$.
Let $X$ be a generically split projective homogeneous $G$-variety
(cf. Example~\ref{intmot}).
By Theorem~\ref{thmMain} we obtain the decomposition
\begin{equation}\label{mndec}
\M(X;\zz/2)\simeq \bigoplus_{i\ge 0} \RR_p(G)(i)^{\oplus a_i},
\end{equation}
where the motive $\RR_p(G)$ is indecomposable and 
(cf. \cite[(5.4-5.5)]{Vo03})
$$
\overline{\RR_p(G)}\simeq 
\bigoplus_{i=0}^{p-1} (\zz/p)(i\cdot (p+1)).
$$

We now identify $\RR_p(G)$.
Let $\rr$ be the Rost invariant as defined
in \cite{Me03} and $\rr_p$ denote its restriction
to the $p$-primary closure of $F$.

\begin{lem}\label{l74}
Let $G$ be a simple linear algebraic group over $F$
satisfying $r=1$ and $d_1>1$ and $p$ be a prime. 
Then $\rr_p(G)$ is trivial iff $\RR_p(G)\simeq\zz/p$.
\end{lem}
\begin{proof}
According to \cite[Theorem~0.5]{Ga01}, \cite{Ch94} and 
\cite[Theoreme~10]{Gi00}
the invariant $\rr_p(G)$ is trivial iff 
the group $G$ splits over the $p$-primary closure of $F$. 
By Corollary~\ref{trivJ}
the latter is equivalent to the fact that $\RR_p(G)\simeq\zz/p$.
\end{proof}

\begin{lem}\label{roinv} Let $G$ and $G'$ be simple linear algebraic
groups over $F$
satisfying $r=1$ and $d_1>1$.
If $\rr_p(G)=\rr_p(G')c$ for some $c\in(\zz/p)^\times$, 
then $\RR_p(G)\simeq \RR_p(G')$.
\end{lem}
\begin{proof} By transfer arguments (see Proposition~\ref{thmComp})
it is enough to prove this over a $p$-primary closure of $F$. 
Let $X$ and $X'$ be the respective varieties of complete flags.
Observe that the invariant $\rr_p(G)$ becomes trivial over 
the function field $F(X)$. Since $\rr_p(G)=\rr_p(G')c$, it becomes
trivial over $F(X')$ as well. 
By Lemma~\ref{l74} $X$ splits over $F(X')$. Similarly $X'$
splits over $F(X)$.

Therefore by Lemma~\ref{gensp}
there exists a cycle $\phi$ in 
$\overline{\Ch_{\dim X}}(X\times X')$ of the form
$\phi=1\times pt + \sum_{\codim \alpha_i>0} 
\alpha_i\times\beta_i$. 
Observe that by definition
$\phi_\star\colon pt_X\mapsto pt_{X'}$.
Similarly, interchanging $X$ and $X'$ 
we obtain a cycle $\phi'\in \overline{\Ch_{\dim X'}}(X'\times X)$ 
such that $\phi_\star'\colon pt_{X'}\mapsto pt_X$.
Restricting of $\phi$ and $\phi'$ to the direct summands
$\overline{\RR_p(G)}$ and $\overline{\RR_p(G')}$ 
of $\M(\BX)$ and $\M(\BX')$ respectively 
we obtain the rational maps 
$\phi_R\colon \overline{\RR_p(G)}\to\overline{\RR_p(G')}$
and $\phi_R'\colon \overline{\RR_p(G')}\to\overline{\RR_p(G)}$.

Since the motive $\RR_p(G)$ is indecomposable and 
$\rk \Ch^i(\overline{\RR_p(G)})\le 1$ for all $i$,
the ring of rational endomorphisms of $\overline{\RR_p(G)}$ is generated
by the identity endomorphism $\Delta$.
The same holds for the ring of rational endomorphisms of 
$\overline{\RR_p(G')}$.
Since $(\phi_R')_\star\circ (\phi_R)_\star\colon pt_X\mapsto pt_X$,
the composite $\phi_R'\circ\phi_R=\Delta$.
Similarly we obtain $\phi_R\circ\phi_R'=\Delta'$.
By Rost Nilpotence, since $\phi_R$ and $\phi_R'$ are rational, the motives
$\RR_p(G)$ and $\RR_p(G')$ are isomorphic.
\end{proof}

\paragraph{\it $\zz$-coefficients.}
Let $G$ be a group of type $\F$ or a strongly inner form of type $\E_6$
which doesn't split 
by field extensions of degrees $2$ and $3$. Observe that such
a group splits by an extension of degree $6$.
Let $X$ be a generically split projective homogeneous $G$-variety.
Then according to Proposition~\ref{intd}
the Chow motive of $X$ with integer coefficients 
splits as a direct sum
of twisted copies of an indecomposable motive $\RR(G)$ such that
\begin{align*}
\RR(G)\otimes \zz/2&=\bigoplus_{i=0,1,2,6,7,8} \RR_2(G)(i),\qquad 
P(\overline{\RR_2(G)},t)=1+t^3,
\\
\RR(G)\otimes \zz/3&=\bigoplus_{i=0,1,2,3}\RR_3(G)(i),\qquad
P(\overline{\RR_3(G)},t)=1+t^4+t^8,
\\
P(\overline{\RR(G)},t)&=1+t+t^2+\ldots+t^{11}.
\end{align*}

\begin{rem} In particular, we provided a uniform
proof of the main results of papers \cite{Bo03} and \cite{NSZ}, where
the cases of $\G$- and $\F$-varieties were considered.
\end{rem}

\begin{rem}
Using Proposition~\ref{intd} one can construct other liftings
of the motivic decompositions of $X$. Thus, the Krull-Schmidt theorem
fails in the category of Chow motives with $\zz/6$-coefficients.
\end{rem}

\paragraph{The case $r>1$.}
According to Table~\ref{bigtab} this holds for groups $G$ of
types $\B_n$ and $\D_n$ and exceptional types 
$\E_7$, $\E_8$ for $p=2$ and
$\E_6^{ad}$, $\E_8$ for $p=3$.

\paragraph{\it Pfister case.}
Let $G=\OO(\phi)$, where $\phi$ is a $k$-fold Pfister form or 
its maximal neighbor. 
Assume $J_2(G)\ne(0,\ldots,0)$. In view of Corollary~\ref{trivJ}
this holds iff $n_2(G)\neq 1$. By Springer's Theorem the latter holds
iff $\phi$ is not split.
By Theorem~\ref{thmMain} we obtain the decomposition
$$
\M(X;\zz/2)\simeq \bigoplus_{i\ge 0} \RR_2(G)(i)^{\oplus a_i}
$$
where $\RR_2(G)$ is indecomposable. Moreover,
by Theorem~\ref{lem47} the same decomposition holds
with $\zz$-coefficients. 

Now we compute $J_2(G)$. 
Let $Y$ be a projective quadric corresponding to $\phi$. 
Then $G$ splits over $F(Y)$ and $Y$ splits over $F(x)$ for any $x\in X$. 
It is known that $\CH^l(\overline{Y})$ for $l<2^{k-1}-1$ 
is generated by $\CH^1(\overline{Y})$ and, therefore, is rational. 
By Proposition~\ref{thmKarp} and Table~\ref{bigtab}
we see that $j_i(G)=0$ for $0\le i<r$, where $r=2^{k-2}$. 
Therefore, $J_2(G)=(0,\ldots,0,1)$ and 
$P(\overline{\RR_2(G)},t)=1+t^{2^{k-1}-1}$.
Finally, by Corollary~\ref{remKrull} the motive $\RR_2(G)$
coincides with the motive introduced in \cite{Ro98} 
which is called the {\em Rost motive}. 

In this way we obtain the Rost decomposition
of the motive of a Pfister quadric and of its maximal neighbor.

\paragraph{\it Maximal orthogonal Grassmannian.}
Let $G=\OO(q)$, where 
$q\colon V\to F$ 
is an arbitrary anisotropic regular quadratic form and $X$ is
a connected component of 
the respective maximal orthogonal Grassmannian. 
The variety $X$ is generically split, hence, 
by Theorem~\ref{thmMain} we have the decomposition
$$
\M(X;\zz/2)\simeq \bigoplus_{i\ge 0} \RR_2(G)(i)^{\oplus a_i},
$$
where the motive $\RR_2(G)$ is indecomposable. 
Comparing the Poincar\'e polynomials of $\M(X;\zz/2)$ and $\RR_2(G)$ we 
obtain the following particular cases:
\begin{itemize} 
\item If the group $G$ corresponds to
a generic cocycle (see \ref{notgen}), 
the motive $\M(X;\zz/2)$ is isomorphic to $\RR_2(G)$ and, hence, 
is indecomposable. This corresponds to the maximal value of the $J$-invariant.
\item If $q$ is a Pfister form or its maximal neighbor, by the previous example
$\RR_2(G)$ coincides
with the Rost motive. This corresponds to the minimal
non-trivial value of the $J$-invariant.
\end{itemize}


\bibliographystyle{chicago}

\noindent
{\sc Viktor Petrov\\
PIMS, Department of Mathematical and Statistical Sciences, University of Alberta,
Edmonton, AB T6G 2G1, Canada}

\ 

\noindent
{\sc Nikita Semenov\\
Mathematisches Institut der LMU M\"unchen, Germany}

\

\noindent
{\sc Kirill Zainoulline\\
Mathematisches Institut der LMU M\"unchen,\\ 
Theresienstr.~39, D-80333 M\"unchen, Germany}

\end{document}